%%%%%%%%%%%%%%%%%%%%%%%%%%%%%%%%%%%%%%%%%%%%%%%%%%%%%%%%%%
%%%%%%%%%%%%%%%%%%%%%%%%%%%%%%%%%%%%%%%%%%%%%%%%%%%%%%%%%%
%%
%%     This is the AMS-LaTeX file:
%%
%%     Colli-Gilardi-Sprekels 25
%%     A distributed control problem for a fractional tumor growth model
%%     
%%
%%%%%%%%%%%%%%%%%%%%%%%%%%%%%%%%%%%%%%%%%%%%%%%%%%%%%%%%%%

\def\LaTeX{%
  \let\Begin\begin
  \let\End\end
  \let\salta\relax
  \let\finqui\relax
  \let\futuro\relax}

\def\UK{\def\our{our}\let\sz s}
\def\USA{\def\our{or}\let\sz z}

\UK
%\USA

%%%%%%%%%%%%%%%%%%%%%%%%%%%%%%%%%

% scegliere fra \TeX e \LaTeX  e fra   uk oppure \USA

%\TeX
\LaTeX

% uk
\USA

%%%%%%%%%%%%%%%%%%%%%%%%%%%%%%%%%
%% page layout
%%%%%%%%%%%%%%%%%%%%%%%%%%%%%%%%%

\salta

\documentclass[twoside,12pt]{article}
\setlength{\textheight}{24cm}
\setlength{\textwidth}{16cm}
\setlength{\oddsidemargin}{2mm}
\setlength{\evensidemargin}{2mm}
\setlength{\topmargin}{-15mm}
\parskip2mm

%%%%%%%%%%%%%%%%%%%%%%%%%%%%%%%%%
%% packages
%%%%%%%%%%%%%%%%%%%%%%%%%%%%%%%%%

%\usepackage{color}
\usepackage[usenames,dvipsnames]{color}
\usepackage{amsmath}
\usepackage{amsthm}
\usepackage{amssymb}
\usepackage[mathcal]{euscript}

%%%%%%%%%%%%%%%%%%%%%%%%%%%%%%%%%
%% added by Pier
%%%%%%%%%%%%%%%%%%%%%%%%%%%%%%%%%

\usepackage{cite}

%%%%%%%%%%%%%%%%%%%%%%%%%%%%%%%%%

%\usepackage[notref,notcite]{showkeys}
%\usepackage{showkeys}
%
%		COLORS FOR CORRECTIONS
%
% do the same, please (i.e., don't use the standard {\color{red} text} or similar): 
% just choose the color you prefer in \def\yourname

% example of use:  \juerg{I want this to become blue}

\definecolor{viola}{rgb}{0.3,0,0.7}
\definecolor{ciclamino}{rgb}{0.5,0,0.5}

\def\pier #1{#1}
\def\juerg #1{#1}
\def\gianni #1{#1}

%%%%%%%%%%%%%%%%%%%%%%%%%%%%%%%%%
%% you may adjust the baseline
%%%%%%%%%%%%%%%%%%%%%%%%%%%%%%%%%

%\renewcommand{\baselinestretch}{0.95}

%%%%%%%%%%%%%%%%%%%%%%%%%%%%%%%%%
%% bibliographystyle
%%%%%%%%%%%%%%%%%%%%%%%%%%%%%%%%%

\bibliographystyle{plain}

%%%%%%%%%%%%%%%%%%%%%%%%%%%%%%%%%
%% environments
%%%%%%%%%%%%%%%%%%%%%%%%%%%%%%%%%

%

\finqui

\def\Beq{\Begin{equation}}
\def\Eeq{\End{equation}}

\def\Bthm{\Begin{theorem}}
\def\Ethm{\End{theorem}}
\def\Blem{\Begin{lemma}}
\def\Elem{\End{lemma}}
\def\Bprop{\Begin{proposition}}
\def\Eprop{\End{proposition}}

\def\Brem{\Begin{remark}\rm}
\def\Erem{\End{remark}}

\def\Bdim{\Begin{proof}}
\def\Edim{\End{proof}}
\def\Bcenter{\Begin{center}}
\def\Ecenter{\End{center}}
\let\non\nonumber

%%%%%%%%%%%%%%%%%%%%%%%%%%%%%%%%%
%% macros
%%%%%%%%%%%%%%%%%%%%%%%%%%%%%%%%%

% macro salvate

% sottosezioni non numerate

\def\step #1 \par{\medskip\noindent{\bf #1.}\quad}

% abbreviazioni di parole

\def\Lip{Lip\-schitz}
\def\Holder{H\"older}
\def\Frechet{Fr\'echet}
\def\aand{\quad\hbox{and}\quad}

\def\lhs{left-hand side}
\def\rhs{right-hand side}

% versioni inglesi (UK) o americane (USA)

% bold, cal e mathop

\def\multibold #1{\def\arg{#1}%
  \ifx\arg\pto \let\next\relax
  \else
  \def\next{\expandafter
    \def\csname #1#1#1\endcsname{{\bf #1}}%
    \multibold}%
  \fi \next}

\def\pto{.}

\def\multical #1{\def\arg{#1}%
  \ifx\arg\pto \let\next\relax
  \else
  \def\next{\expandafter
    \def\csname cal#1\endcsname{{\cal #1}}%
    \multical}%
  \fi \next}

% operatori

\def\multimathop #1 {\def\arg{#1}%
  \ifx\arg\pto \let\next\relax
  \else
  \def\next{\expandafter
    \def\csname #1\endcsname{\mathop{\rm #1}\nolimits}%
    \multimathop}%
  \fi \next}

\multibold
qwertyuiopasdfghjklzxcvbnmQWERTYUIOPASDFGHJKLZXCVBNM.

\multical
QWERTYUIOPASDFGHJKLZXCVBNM.

\multimathop
diag dist div dom mean meas sign supp .

% accorpamenti di formule citate:
% uso  \accorpa {prima}{seconda}
%      \Accorpa\cs prima seconda (con il comodo blank anche dopo)
% NB: \Accorpa definisce \cs come l'accorpamento delle due citazioni
% e scrive sul file.log

\def\accorpa #1#2{\eqref{#1}--\eqref{#2}}
\def\Accorpa #1#2 #3 {\gdef #1{\eqref{#2}--\eqref{#3}}%
  \wlog{}\wlog{\string #1 -> #2 - #3}\wlog{}}

% macro comode

\def\separa{\noalign{\allowbreak}}

\def\somma #1#2#3{\sum_{#1=#2}^{#3}}

\def\graffe #1{\mathopen\{#1\mathclose\}}

\def\<#1>{\mathopen\langle #1\mathclose\rangle}
\def\norma #1{\mathopen \| #1\mathclose \|}

\def\[#1]{\mathopen\langle\!\langle #1\mathclose\rangle\!\rangle}

\def\iot {\int_0^t}
\def\ioT {\int_0^T}
\def\intQt{\int_{Q_t}}
\def\intQ{\int_Q}
\def\iO{\int_\Omega}
\def\itT{\int_t^T}
\def\bintQt{\itT\!\!\!\iO}

\def\dt{\partial_t}
\def\dn{\partial_\nu}

\def\cpto{\,\cdot\,}

\def\checkmmode #1{\relax\ifmmode\hbox{#1}\else{#1}\fi}
\def\aeO{\checkmmode{a.e.\ in~$\Omega$}}
\def\aeQ{\checkmmode{a.e.\ in~$Q$}}

\def\aet{\checkmmode{a.e.\ in~$(0,T)$}}

\def\aat{\checkmmode{for a.a.~$t\in(0,T)$}}

% insiemi numerici

\def\erre{{\mathbb{R}}}

\def\enne{{\mathbb{N}}}

% spazi di funzioni a valori vettoriali su [0,T], [0,t], [0,s], [0,+\infty), [\delta,T]

% Come ricordare: in generale i simboli L H W  C da soli per gli spazi su (0,T)
% gli stessi raddoppiati per (0,+\infty)
% aggiunta di t o s al simbolo per (0,t) e (0,s)
% aggiunta di d al simbolo semplice o doppio per intervalli (\delta,T) e (\delta,+\infty)
% il simbolo C e i suoi derivati mettono le quadre anziche' le tonde

% Esempi   \L2V   \L\infty\Vp   \W{1,1}H   \C0H   \LL2V   \CC0\Vp   \Ld2V  \CCdH

\def\genspazio #1#2#3#4#5{#1^{#2}(#5,#4;#3)}
\def\spazio #1#2#3{\genspazio {#1}{#2}{#3}T0}

\def\L {\spazio L}
\def\H {\spazio H}
\def\W {\spazio W}

\def\C #1#2{C^{#1}([0,T];#2)}

% spazi di funzioni su \Omega, \Gamma, Q e \Sigma

\def\Lx #1{L^{#1}(\Omega)}
\def\Hx #1{H^{#1}(\Omega)}

\def\Luno{\Lx 1}
\def\Ldue{\Lx 2}
\def\Linfty{\Lx\infty}

\def\Huno{\Hx 1}
\def\Hdue{\Hx 2}
\def\Hunoz{{H^1_0(\Omega)}}

% spazi di funzioni su Q e S

\def\LQ #1{L^{#1}(Q)}

% lettere greche

%\let\theta\vartheta

\let\phi\varphi

\let\TeXchi\chi                         % new \chi, exactly on the baseline
\newbox\chibox
\setbox0 \hbox{\mathsurround0pt $\TeXchi$}
\setbox\chibox \hbox{\raise\dp0 \box 0 }
\def\chi{\copy\chibox}

% quadratino di fine dimostrazione

% abbreviazioni specifiche del lavoro

\def\VA #1{V_A^{#1}}
\def\VB #1{V_B^{#1}}
\def\VC #1{V_C^{#1}}

\def\Ar{{A,\,\rho}}
\def\Bs{{B,\,\sigma}}
\def\Ct{{C,\,\tau}}

\def\phih{\hat\phi_h}
\def\muh{\hat\mu_h}
\def\Sh{\hat S_h}

\def\phiz{\phi_0}
\def\Sz{S_0}

\def\soluz{(\mu,\phi,S)}
\def\soluzb{(\mub,\phib,\Sb)}
\def\soluzl{(\eta,\xi,\zeta)}
\def\soluza{(q,p,r)}
\def\ds{\,ds}
\def\dx{\,dx}

\def\Uad{\calU_{ad}}
\def\umin{u_{min}}
\def\umax{u_{max}}
\def\BR{\calB_R}
\def\az{a_0}
\def\bz{b_0}
\def\aM{a_M}
\def\bM{b_M}
\def\aR{a_R}
\def\bR{b_R}

\def\ub{\overline u}
\def\mub{\overline\mu}
\def\phib{\overline\phi}
\def\Sb{\overline S}
\def\muh{\mu^h}
\def\phih{\phi^h}
\def\Sh{S^h}
\def\xih{\xi^h}
\def\etah{\eta^h}
\def\zetah{\zeta^h}
\def\Qh #1{Q^h_#1}
\def\Rh #1{R^h_#1}
\def\un{u_n}
\def\mun{\mu_n}
\def\phin{\phi_n}
\def\Sn{S_n}
\def\qn{q^n}
\def\pn{p^n}
\def\rn{r^n}
\def\soluzn{(\qn,\pn,\rn)}

\def\VAn{\VA{(n)}}
\def\VBn{\VB{(n)}}
\def\VCn{\VC{(n)}}
\def\VAm{\VA{(m)}}
\def\VBm{\VB{(m)}}

%%%%%%%%%%%%%%%%%%%%%%%%%%%%%%
\Begin{document}
%%%%%%%%%%%%%%%%%%%%%%%%%%%%%%%%%

%%%%%%%%%%%%%%%%%%%%%%%%%%%%%%%%%
%% front page
%%%%%%%%%%%%%%%%%%%%%%%%%%%%%%%%%

%
\title{A distributed control problem\\for a fractional tumor growth model}
\author{}
\date{}
\maketitle
\Bcenter
\vskip-1cm
{\large\sc Pierluigi Colli$^{(1)}$}\\
{\normalsize e-mail: {\tt pierluigi.colli@unipv.it}}\\[.25cm]
{\large\sc Gianni Gilardi$^{(1)}$}\\
{\normalsize e-mail: {\tt gianni.gilardi@unipv.it}}\\[.25cm]
{\large\sc J\"urgen Sprekels$^{(2)}$}\\
{\normalsize e-mail: {\tt sprekels@wias-berlin.de}}\\[.45cm]
$^{(1)}$
{\small Dipartimento di Matematica ``F. Casorati'', Universit\`a di Pavia}\\
{\small and Research Associate at the IMATI -- C.N.R. Pavia}\\
{\small via Ferrata 5, 27100 Pavia, Italy}\\[.2cm]
$^{(2)}$
{\small Department of Mathematics}\\
{\small Humboldt-Universit\"at zu Berlin}\\
{\small Unter den Linden 6, 10099 Berlin, Germany}\\[2mm]
{\small and}\\[2mm]
{\small Weierstrass Institute for Applied Analysis and Stochastics}\\
{\small Mohrenstrasse 39, 10117 Berlin, Germany}
\Ecenter
\Begin{abstract}\noindent
\juerg{In this paper, we study the distributed optimal control of a system of three
evolutionary equations involving fractional powers of three selfadjoint, monotone,
unbounded linear operators having compact resolvents. The system is a generalization of a 
Cahn--Hilliard type phase field system modeling tumor growth that goes back to Hawkins-Daarud et~al.\ 
({\em Int.\ J.\ Numer.\ Math.\ Biomed.\ Eng.\ {\bf 28} (2012), 3--24.}) The aim of the control process,
which could be realized by either administering a drug or monitoring the nutrition, is to keep the
tumor cell fraction under control while avoiding possible harm for the patient. In contrast to
previous studies, in which the occurring unbounded operators governing the diffusional regimes
were all given by the Laplacian with zero Neumann boundary conditions, the operators may in our case be different; more generally, we consider systems with  fractional powers of the type that were studied
in the recent work {\em Adv.\ Math.\ Sci.\ Appl.\ {\bf 28} (2019), 343--375,} by the present authors.          
In our analysis, we show the Fr\'echet differentiability of the associated control-to-state 
operator, establish the existence of solutions to the associated adjoint system, and derive
the first-order necessary conditions of optimality for a cost functional of tracking type. }

\vskip3mm
\noindent {\bf Key words:}
Fractional operators, Cahn--Hilliard systems, well-posedness, regularity,
\juerg{optimal control, necessary optimality conditions.}
\vskip3mm
\noindent {\bf AMS (MOS) Subject Classification:} \juerg{35K55, 35Q92, 49J20, 92C50.}
\End{abstract}
\salta
\pagestyle{myheadings}
\newcommand\testopari{\pier{\sc Colli \ --- \ Gilardi \ --- \ Sprekels}}
\newcommand\testodispari{\pier{\sc Optimal control of a fractional tumor growth model}}
\markboth{\testopari}{\testodispari}
\finqui
%
%%%%%%%%%%%%%%%%%%%%%%%%%%%%%%%%%
%% very beginning
%%%%%%%%%%%%%%%%%%%%%%%%%%%%%%%%%

\section{Introduction}
\label{Intro}
\setcounter{equation}{0}

The recent paper \cite{CGS23} investigates the evolutionary system
\begin{align}
  & \alpha \,\dt\mu + \dt\phi + A^{2\rho} \mu 
  = P(\phi)(S-\mu),
  \label{Iprima}
  \\
  & \beta \,\dt\phi + B^{2\sigma}\phi + f(\phi) = \mu,
  \label{Iseconda}
  \\
  & \dt S + C^{2\tau}S 
  = - P(\phi)(S-\mu) + u,
  \label{Iterza}
\end{align}
where the equations are understood to hold in $\Omega$,
a bounded, connected and smooth domain in~$\erre^3$,
and in the time interval~$(0,T)$.
In the above system, $A^{\pier{2\rho}}$, $B^{2\sigma}$, and~$C^{\pier{2\tau}}$, with $r>0$, $\sigma>0$, $\rho>0$,
denote fractional powers of the selfadjoint,
monotone, and unbounded, linear operators $A$, $B$ and~$C$, respectively, 
which are supposed to be densely defined in $H:=\Ldue$ and to have compact resolvents.
Moreover, $\alpha$ and $\beta$ are positive real parameters. 

\juerg{The system \eqref{Iprima}--\eqref{Iterza} is a generalization of a diffuse interface model for tumor growth. 
Such models, which are usually established in the framework of the Cahn--Hilliard model originating from the theory of
phase transitions, have drawn increasing attention in the past years among mathematicians and applied scientists.
We cite here just \cite{CLLW,CL,HDPZO,HDZO,OHP,WLFC,WZZ} as a sample of pioneering papers in this direction.
In this connection, $\,\phi\,$ stands for an order parameter that should attain its values in the interval $[-1,1]$, where
the values $-1$ and $+1$ indicate the healthy cell and tumor cell cases, respectively. The variable $S$ \pier{represents}
the nutrient extra-cellular water concentration, $u$ stands for a source term that acts as a control to monitor the
evolution of the tumor cell fraction $\phi$, and the nonlinearity $P$ occurring in \eqref{Iprima} and \eqref{Iterza}
is a nonnegative and smooth function modeling a proliferation rate. Finally,
$\mu$ represents the chemical potential, which acts as the driving thermodynamic force of the evolution 
and  is obtained as the  variational derivative with respect to the order parameter $\phi$ of a suitable
 free energy functional. In this connection, the nonlinearity $f$~denotes the derivative of a 
double-well potential~$F$ which plays the role of a specific local free energy and yields the main contribution to the total
free energy.}
Important examples for $\,F\,$ are the so-called {\em classical regular potential} and the {\em logarithmic double-well potential\/},
given by the formulas
\begin{align}
  & F_{reg}(r) := \frac 14 \, (r^2-1)^2 \,,
  \quad r \in \erre, 
  \quad \aand
  \label{regpot}
  \\
  & F_{log}(r) := \bigl( (1+r)\ln (1+r)+(1-r)\ln (1-r) \bigr) - c_1 r^2 \,,
  \quad r \in (-1,1),
  \label{logpot}
\end{align}
respectively.
In \eqref{logpot}, the constant $c_1$ is larger than~$1$, so that $F_{log}$ is nonconvex.
Furthermore, the function $P$ in \eqref{Iprima} and \eqref{Iterza}
is nonnegative and smooth.
Finally, the \pier{datum} $u$ appearing in \eqref{Iterza} is given.

\juerg{In the literature, the diffusional developments in the system have usually been modeled by the Laplacian,
that is, the case $A^{2\rho}=B^{2\sigma}=C^{2\tau}=-\Delta$, accompanied by zero Neumann boundary conditions, 
was assumed, where two main classes of models were considered. The first class of models regards the tumor and
healthy cells as inertialess fluids; in such models special fluid effects can be incorporated by 
postulating a Darcy or Stokes--Brinkman law, see, e.g., the works \cite{DFRGM,EGAR,FLRS,GARL_3,GARL_1,
GARL_2, GARL_4, GAR,GLSS,SW, WLFC}, where we also refer to \cite{ConGio,DGG}. The other class of models, to
which the model considered here belongs, neglects the velocity. Typical contributions in this direction were
given in \cite{CRW,CGH,CGRS_VAN,CGRS_OPT,CGRS_ASY,FGR}, to name just a few.}

\juerg{
While the occurrence of more general diffusional regimes of fractional type has been studied for a long time
in the mathematical literature, it was only recently (see, e.g., \cite{AM,ASS, 
CGS19,CGS21,CGS18,CGS22,Gal1,Gal2,Gal3}) that fractional operators 
have been investigated in the framework of Cahn--Hilliard systems (for phase field systems of Caginalp type, see also
\cite{CG}), and the only investigations 
of tumor growth models involving fractional diffusive regimes such as in the system \eqref{Iprima}--\eqref{Iterza} 
seem to be the recent papers \cite{CGS23,CGS24} by the present authors. }

\juerg{
In particular, in the paper \cite{CGS23}, 
under rather general assumptions on the operators and the potentials,}
well-posedness and regularity results for the initial value problem for \accorpa{Iprima}{Iterza} 
were established in the case $\,u=0\pier{,}\,$ \juerg{under the assumption that} $\,\alpha>0\,$ and $\,\beta>0$.
However, some remarks on more general cases including $u\in\LQ2$, where $Q:=\Omega\times(0,T)$, 
have been given \pier{in \cite{CGS23}.}
In particular, under suitable assumptions on the initial data,
for every $u\in\LQ2$, there exists at least a solution $\soluz$ in~a proper functional space
to~a weak version of \accorpa{Iprima}{Iterza}
(namely, \eqref{Iseconda} is replaced by a variational inequality involving 
the convex part $F_1$ of $F$ rather than~$f$, 
since $F_1$ is not supposed to be differentiable).
Moreover, the solution is unique if the domains of the fractional operators $A^\rho$ and $C^\tau$ 
satisfy suitable embeddings of Sobolev type.
Finally, if $B^{2\sigma}$ behaves like the Laplace operator with either Dirichlet or Neumann zero boundary conditions 
and $f$ is single valued (like \eqref{regpot} and~\eqref{logpot}), \juerg{then}
the solution solves equation \eqref{Iseconda} in a stronger sense,
and it is even smoother under more restrictive assumptions on the initial data.

In this paper, we first establish similar results for system \accorpa{Iprima}{Iterza}
by assuming $\alpha=0$ and $\beta>0$ (in~fact, we take $\beta=1$ without loss of generality). \juerg{In
particular, we extend \pier{some} results shown for this case in the recent paper \cite{CGS24}.} 
Then, we discuss a distributed control problem for the modified system.
Namely, given nonnegative constants $\kappa_i$, $i=1,\dots,5$, 
and functions $\phi_Q,S_Q\in\LQ2$ and $\phi_\Omega,S_\Omega\in\Ldue$, 
we consider the problem of minimizing the cost functional
\begin{align*}
  \calJ(u,\phi,S)
  & := \,\juerg{\frac{\kappa_1}2} \intQ |\phi-\phi_Q|^2
  + \juerg{\frac{\kappa_2}2} \iO |\phi(T) - \phi_\Omega|^2
  \non
  \\
  & \hspace*{5mm} + \,\juerg{\frac{\kappa_3}2} \intQ |S-S_Q|^2
  + \juerg{\frac{\kappa_4}2} \iO |S(T) - S_\Omega|^2
  + \juerg{\frac{\kappa_5}2} \intQ |u|^2\,,
  \end{align*}
where $\phi$ and $S$ are the components of the solution $\soluz$ corresponding to the control~$u$, 
which is supposed to vary under restrictions of the type
$\,\umin\leq u\leq \umax$.

\juerg{The choice of this tracking-type cost functional reflects the plan of a medical treatment via the 
application of drugs over some finite time interval $(0,T)$ with the aim of monitoring the evolution of the
tumor fraction $\phi$ under the restriction that no harm be inflicted on the patient. We remark at this 
place that it would be desirable to minimize the duration, i.e., the time $T>0$, of the medical treatment as well, 
in order to prevent that the tumor cells develop a resistance against the drug. However, such an approach, which was possible
(see, e.g.,  \cite{CRW}) in the special case when $A^{2\rho}=B^{2\sigma}=C^{2\tau}=-\Delta$, becomes very
complicated in the situation considered here and was therefore not included.}  

\juerg{The literature on optimal control problems for Cahn--Hilliard systems is still scarce. In this connection, we refer
 the reader to \cite{CGS14,CGS19}, where a number of references \pier{is} given. Even less investigations have been made
  on  optimal control problems for tumor growth models. \pier{About that, let us} refer to the works\pier{\cite{CRW,SM,CGRS_OPT,EK_ADV,EK,GARLR,S_a,S_b,S_DQ,S,SW}}, for various models involving the Laplacian. Concerning the optimal
control of Cahn--Hilliard systems with fractional operators, we just can cite \cite{CGS19,CGS21}, and, to the \gianni{authors'}
best knowledge, the present paper is the first contribution on the optimal control of the tumor growth model with 
fractional operators.}

The remainder of the paper is organized as follows. 
In the next section, we list our assumptions and notations and 
present our results on the state system.
The next Section~\ref{FRECHET} is devoted to the study of the control-to-state mapping
and \pier{of its \Frechet\ differentiability}.
In the last section we deal with the control problem.
Namely, the existence of an optimal control is proved
and the first order necessary conditions involving a proper adjoint system are derived.

%%%%%%%%%%%%%%%%%%%%%%%%%%%%%%%%%%%%%%%%%%%%%%%%%%%%%%%%%%%%%%%%%%%%%%%%

\section{The state system}
\label{STATESYSTEM}
\setcounter{equation}{0}

In this section, we first introduce the notations and the assumptions needed for the analysis of the state system. Then, we present our results.
We closely follow~\cite{CGS23}.
First of all, the set $\Omega\subset\erre^3$ is assumed to be bounded, connected and smooth,
with volume $|\Omega|$ and outward unit normal vector field $\,\nu\,$ on $\Gamma:=\partial\Omega$.
Moreover, $\dn$ stands for the corresponding normal derivative.
We set
\Beq
  H := \Ldue
  \label{defH}
\Eeq
and denote by $\norma\cpto$ and $(\cpto,\cpto)$ the standard norm and inner product of~$H$.
As for the operators, we first postulate that
\begin{align}
  & A:D(A)\subset H\to H , \quad
  B:D(B)\subset H\to H
  \aand
  C:D(C)\subset H\to H
  \quad \hbox{are}
  \non
  \\
  & \hbox{unbounded monotone selfadjoint linear operators with compact resolvents.} 
  \qquad
  \label{hpABC} 
\end{align}
Therefore, there are sequences 
$\{\lambda_j\}$, $\{\lambda'_j\}$, $\{\lambda''_j\}$ \,and\, $\{e_j\}$, $\{e'_j\}$, $\{e''_j\}$ 
\,of eigenvalues and of corresponding eigenvectors satisfying
\begin{align}
  & A e_j = \lambda_j e_j, \quad
  B e'_j = \lambda'_j e'_j,
  \aand
  C e''_j = \lambda''_j e''_j,
  \non
  \\
  & \quad \hbox{with} \quad
  (e_i,e_j) = (e'_i,e'_j) = (e''_i,e''_j) = \delta_{ij}
  \quad \hbox{for $i,j=1,2,\dots$},
  \label{eigen}
  \\
  \separa
  & 0 \leq \lambda_1 \leq \lambda_2 \leq \dots , \quad
  0 \leq \lambda'_1 \leq \lambda'_2 \leq \dots
  \aand
  0 \leq \lambda''_1 \leq \lambda''_2 \leq \dots,  
  \non
  \\
  & \quad \hbox{with} \quad
  \lim_{j\to\infty} \lambda_j
  = \lim_{j\to\infty} \lambda'_j
  = \lim_{j\to\infty} \lambda''_j
  = + \infty,
  \label{eigenvalues}
  \\[1mm]
  & \hbox{$\{e_j\}$, $\{e'_j\}$ \,and\, $\{e''_j\}$ \,are complete systems in $H$}.
  \label{complete}
\end{align}
As a consequence, we can define the powers of the above operators with arbitrary positive real exponents.
As far as the first one is concerned, we have, for $\rho>0$,
\begin{align}
  & \VA\rho := D(A^\rho)
  = \Bigl\{ v\in H:\ \somma j1\infty |\lambda_j^\rho (v,e_j)|^2 < +\infty \Bigr\}
  \aand
  \label{defdomAr}
  \\[-3mm]
  & A^\rho v = \somma j1\infty \lambda_j^\rho (v,e_j) e_j
  \quad \hbox{for $v\in\VA\rho$}\,,
  \label{defAr}
\end{align}
and we endow $\VA\rho$ with the graph norm
\Beq
  \norma v_\Ar := \bigl( \norma v^2 + \norma{A^\rho v}^2 \bigr)^{1/2}
  \quad \hbox{for every $v\in\VA\rho$}.
  \label{defnormaAr}
\Eeq
Similarly, we set
\Beq
  \pier{\VB\sigma := D(B^\sigma)
  \aand
  \VC\tau := D(C^\tau),}
  \label{defBsCt}
\Eeq
\pier{with the graph norms}
\begin{align}
  & \norma v_\Bs := \bigl( \norma v^2 + \norma{B^\sigma v}^2 \bigr)^{1/2}
  \aand
  \norma v_\Ct := \bigl( \norma v^2 + \norma{C^\tau v}^2 \bigr)^{1/2},
  \qquad
  \non
  \\
  & \quad \hbox{for $v\in \VB\sigma$ and $v\in\VC\tau$, respectively}.
  \label{defnormeBsCt}
\end{align}
From now on, we assume:
\Beq
  \hbox{$\rho$, $\sigma$ and $\tau$ are fixed positive real numbers}.
  \label{hpexponents}
\Eeq
However, we need the further assumptions we list at once.
It is understood that all of the embeddings below are assumed to be continuous.
\begin{align}
  & \hbox{The first eigenvalue $\lambda_1$ of $A$ is strictly positive}.
  \label{hpeigen}
  \\[0.5mm]
  & \VA{2\rho} \subset \Linfty , \quad
  \VA\rho \subset \Lx 4 , \quad
  \VB\sigma \subset \Lx 4 ,
  \aand
  \VC\tau \pier{{}\subset{}} \Lx 4 .  \qquad
  \label{embeddings}
  \\[0.5mm]
  & \psi(v) \in H
  \aand
  \bigl( B^{2\sigma} v , \psi(v) \bigr) \geq 0,
  \quad \hbox{for every $\,v\in\VB{2\sigma}\,$ 
  and every monotone}
	\non\\
	&\quad\hbox{\pier{and} \Lip\ continuous function \,$\psi:\erre\to\erre$\, vanishing at the origin}.
  \label{hpBs}
\end{align}
Due to the continuus embeddings \eqref{embeddings},
there exists a constant $C_*>0$ such that
\begin{align}
  & \norma v_\infty \leq C_* \juerg{\|v\|_{A,2\rho}} \,, \enskip
  \norma v_4 \leq C_* \norma v_\Ar \,, \enskip
  \norma v_4 \leq C_* \norma v_\Bs\,,
  \enskip \hbox{and} \enskip
  \norma v_4 \leq C_* \norma v_\Ct\,,
    \non
  \\
  & \quad \hbox{for every $v\in\VA{2\rho}$, $v\in\VA\rho$, $v\in\VB\sigma$, and 
	$v\in\VC\tau$, \,\,respectively,}
  \label{embconst}
\end{align}
where, for $p\in[1,+\infty]$, the symbol $\norma\cpto_p$ denotes the norm in $\Lx p$.
The same symbol will \juerg{also be used} for the norm in $\LQ p$ provided that no confusion can arise.

\Brem
\label{Remoperators}
We have to make some comments on \accorpa{hpeigen}{hpBs}.
The first of these assumptions is satisfied if \,$A$\, is, e.g., the Laplace operator $-\Delta$
with zero Dirichlet (or~Robin) boundary conditions,
while the case of zero Neumann boundary conditions is excluded \pier{unless one adds to the Laplace 
operator, e.g., some zero-order term ensuring coerciveness}.
However, it is clear that $\,A\,$ could be a much more general operator.
\pier{By still} considering the Laplace operator with (zero) Dirichlet boundary conditions as~$A$,
we can also discuss the first two embeddings in~\eqref{embeddings}.
By noting that $D(A)=\Hdue\cap\Hunoz$ and $\Omega$ is smooth, 
it results that $D(A^{2\rho})\subset\Hx{4\rho}$ and $D(A^\rho)\subset\Hx{2\rho}$.
Hence, both embeddings hold true if $\,\rho\geq 3/8,$ since $\,\Omega\,$ is three-dimensional.
\pier{Finally,} we make a comment on~\eqref{hpBs}.
Assume, for instance, that $B^{2\sigma}=-\Delta$ with zero Neumann boundary conditions.
Then, $\VB{2\sigma}=\{v\in\Hdue:\ \dn v=0\,\,\mbox{ on }\,\Gamma\}$
and, for every $v\in\VB{2\sigma}$ and $\psi$ as in \eqref{hpBs}, 
we have that $\psi(v)\in\Huno$ (since $v\in\Huno$)~and
\Beq
  \bigl( B^{2\sigma} v, \psi(v) \bigr)
  = \iO (-\Delta v) \, \psi(v)
  = \iO \nabla v \cdot \nabla\psi(v)
  = \iO \psi'(v) |\nabla v|^2
  \geq 0 .
  \non
\Eeq
The same argument works if we take the Dirichlet boundary conditions instead of the Neumann ones,
since the functions $\psi$ \juerg{considered in}
\eqref{hpBs} vanish at the origin.
More generally, $B^{2\sigma}$ can be the principal part of an elliptic operator 
in divergence form with smooth coefficients.
\pier{In particular,} even though some restrictions on $A$, $B$, \pier{and $C$}
\juerg{have to be imposed in order to fulfill the} properties \accorpa{hpeigen}{hpBs},
no relationship between them is needed, and the three operators can be completely independent from each other.
\Erem

\Brem
\label{Firsteigenvalue}
Assumption \eqref{hpeigen} allows us to consider an equivalent norm in~$\VA\rho$.
Indeed, for every $\,v\in\VA\rho\,$ we have that
\Beq
  \norma{A^\rho v}^2 
  = \somma j1\infty |\lambda_j^\rho (v,e_j)|^2
  \geq \lambda_1^{2\rho} \somma j1\infty |(v,e_j)|^2
  = \lambda_1^{2\rho} \norma v^2 .
\Eeq
Hence, since $\lambda_1>0$, we deduce that
\Beq
  \norma v \leq \lambda_1^{-\rho} \norma{A^\rho v}
  \quad \hbox{for every $v\in\VA\rho$},
  \label{forequiv}
\Eeq
so that the function $\,v\mapsto\norma{A^\rho v}\,$ defines a norm in $\VA\rho$ 
that is equivalent to the graph norm~\eqref{defnormaAr}.
\Erem

For the nonlinear functions entering our system, we postulate the following properties:
\begin{align}
  & \hbox{$D(F)$ is an open interval $(a,b)$ of the real line with $0\in(a,b)$}.
  \label{hpDF}
  \\
  & F := D(F) \to \erre
  \quad \hbox{is a $C^3$ function.} 
  \qquad
  \label{hpF}
  \\
  & F(s) \geq C_1 s^2 - C_2
  \aand F''(s) \geq -C_3
  \non
  \\
  & \quad \hbox{for some constants $C_i>0$ and every $s\in D(F)$}.
  \label{hpbelow}
  \\
  & \hbox{$f:=F'$ satisfies} \quad
  \lim_{s\searrow a} f(s) = - \infty 
  \aand
  \lim_{s\nearrow b} f(s) = + \infty \,.
  \label{hpf}
  \\
  & P :\erre \to [0,+\infty) \quad \hbox{is bounded and \Lip\ continuous on $\erre$}
  \non
  \\
  & \quad \hbox{and of class $C^2$ in $D(F)$}.
  \label{hpP}
\end{align}
Clearly, \accorpa{hpDF}{hpf} are fulfilled by the \pier{significant} potentials \eqref{regpot} and \eqref{logpot}.

\Brem
\label{HpCGS23}
The hypotheses \accorpa{hpDF}{hpf} on $F$ ensure that the conditions required in \cite{CGS23}\pier{, i.e., $F =\bigl( \tilde F_1 + \tilde F_2 \bigr)|_{(a,b)}$, where
\begin{align}
  & \tilde F_1 : \erre \to [0,+\infty]
  \quad \hbox{is convex, proper, and l.s.c.,\ with} \quad
  \tilde F_1(0) = 0\pier{,}
  \label{hpFuno23}
  \\
  \separa
  & \tilde F_2 : \erre \to \erre
  \quad \hbox{is of class $C^1$ with a \Lip\ continuous first derivative\pier{,}
}
  \label{hpFdue23}
  \\
  & \tilde F_1^\lambda (s) + \tilde F_2(s) \geq - C_0
  \quad \hbox{for some constant $C_0$ and every $s\in\erre$,
}
  \label{hpbelow23}
\end{align}
are satisfied, as we show at once.
We first split $F$ by defining,} for $s\in(a,b)$,
\begin{align}
  & f_1(s) := \int_0^s (f'(s'))^+ \ds', \quad
  F_1(s) := \int_0^s f_1(s') \ds' \,,
  \non
  \\
  & f_2(s) := F'(0) - \int_0^s (f'(s'))^- \ds'
  \aand
  F_2(s) = F(0) + \int_0^s f_2(s') \ds' \,.
  \non
\end{align} 
Notice that $F_1$ is nonnegative and convex and that $F_1(0)=0$.
If $(a,b)\not=\erre$, we properly extend these functions $F_i$ 
to functions $\tilde F_i$ defined in the whole of~$\erre$. 
One can preserve the mentioned properties of $F_1$, including its lower semicontinuity, by setting
\Beq
  \tilde F_1(a) := \lim_{s\searrow a} F_1(s) , \quad
  \tilde F_1(b) := \lim_{s\nearrow b} F_1(s) 
  \aand
  \tilde F_1(s) := +\infty 
  \quad \hbox{for $s\not\in[a,b]$}.
  \non
\Eeq
Moreover, one can ensure that the derivative of the extension $\tilde F_2$ is \Lip\ continuous,
by noting that $F_2'=f_2$ already is \Lip\ continuous in $(a,b)$
since its derivative $f_2'=-(f')^-$ is bounded by the assumption \eqref{hpbelow} on~$F''$.
The last condition that we have to check is \pier{\eqref{hpbelow23}, 
where $\tilde F_1^\lambda$ is the Moreau--Yosida approximation of~$\tilde F_1$ at the level~$\lambda$.}
We notice that this condition is not equivalent to an inequality of type $F(s)\geq-C_2$,
which follows from~\eqref{hpbelow} and looks rather natural in performing formal a priori estimates.
On the other hand, one can prove that 
the inequality we need is implied by the full quadratic growth condition given in~\eqref{hpbelow}
(see \cite[formula (3.1)]{CGS22} for some explanation).
For this reason, we have postulated the latter.
\Erem

Although some of the results \juerg{to be presented  will  not require
the whole set of hypotheses made so far, 
the statements will be greatly simplified if we do not each time recall
the properties of the involved operators, spaces, and nonlinearities; we therefore make the following 
general assumption:}
\Beq
 \hbox{\it All of the \pier{assumptions} made above on the structure are in force from now on.}
  \label{allassumptions}
\Eeq

As mentioned in the Introduction, we only deal with the case $\alpha=0$ and $\beta>0$
of system \accorpa{Iprima}{Iterza}. 
Clearly, we can take $\beta=1$ without loss of generality.
Hence, the \juerg{Cauchy problem forming the state system under investigation reads as follows:}
\begin{align}
  & \dt\phi + A^{2\rho} \mu 
  = P(\phi)(S-\mu)\,,
  \label{prima}
  \\[0.5mm]
  & \dt\phi + B^{2\sigma}\phi + f(\phi) = \mu\,,
  \label{seconda}
  \\[0.5mm]
  & \dt S + C^{2\tau}S 
  = - P(\phi)(S-\mu) + u\,,
  \label{terza}
  \\[0.5mm]
  & \phi(0) = \phiz
  \aand
  S(0) = \Sz,
  \label{cauchy}
\end{align}
\Accorpa\Pbl prima cauchy
where $\phiz$ and $\Sz$ are prescribed initial data that are supposed to satisfy
\begin{align}
  & \phiz \in \VB{2\sigma} \,, \quad
  \Sz \in \VC\tau\,,
  \aand
  \az \leq \phiz \leq \bz \quad \aeO 
  \non
  \\
  & \quad \hbox{for some compact interval $[\az\,,\bz]\subset(a,b)$}.
  \label{hpz}
\end{align}
In fact, we could solve a weak form of the above problem under milder assumption on the initial data;
however, in order to guarantee a sufficient regularity level of the solution, we need the whole of~\eqref{hpz}.
Given a final time $T\in(0,+\infty)$,
the regularity we can ensure (besides some boundedness to be discussed later~on) 
is~the following:
\begin{align}
  & \mu \in \L\infty{\VA{2\rho}}\,,
  \label{regmu}
  \\[0.5mm]
  & \phi \in \W{1,\infty}H \cap \H1{\VB\sigma} \cap \L\infty{\VB{2\sigma}}\,,
  \label{regphi}
  \\[0.5mm]
  & f(\phi) \in \L\infty H\,,
  \label{regfphi}
  \\[0.5mm]
  & S \in \H1H \cap \L\infty{\VC\tau} \cap \L2{\VC{2\tau}}\,, 
  \label{regS}
\end{align}
\Accorpa\Regsoluz regmu regS
so that equations \accorpa{prima}{terza} are satisfied \aeQ, 
where we recall that
\Beq
  Q := \Omega \times (0,T).
  \label{defQ}
\Eeq
We notice at once that the first embedding in \eqref{embeddings} yields that
\Beq
  \mu \in \L\infty{\VA{2\rho}}
  \quad \hbox{implies that} \quad
  \mu \in \LQ\infty .
  \label{mubdd}
\Eeq
At this point, we are ready to present our results.
To this end, it is convenient to introduce the following variational formulation of \accorpa{prima}{terza}:
\begin{align}
  & \bigl( \dt\phi(t) , v \bigr)
  + \juerg{( A^\rho \mu(t) , A^\rho v )}
  = \bigl( P(\phi \juerg{(t)}) ( S(t) - \mu(t) ) , v \bigr)
  \non
  \\
  & \quad \hbox{for every $v\in\VA\rho$ and \aat}\,,
  \label{primav}
  \\[2mm]
  \separa
  & \bigl( \dt\phi(t) , v \bigr)
  + \bigl( B^\sigma \phi(t) , B^\sigma v \bigr)
  + \bigl( f(\phi(t)) ,  v \bigr)
  = \bigl( \mu(t) , v \bigr)
  \non
  \\
  & \quad \hbox{for every $v\in\VB\sigma$ and \aat}\,,
  \label{secondav}
  \\[2mm]
  \separa
  & \bigl( \dt S(t) , v \bigr)
  + \juerg{( C^\tau S(t) ,C^\tau v )}
  = - \bigl( P(\phi \juerg{(t)}) ( S(t) - \mu(t) ) , v \bigr)\juerg{\,+\,(u(t),v)}
  \non
  \\
  & \quad \hbox{for every $v\in\VC\tau$ and \aat}.
  \label{terzav}
\end{align}
This is based on obvious properties of the powers of the operators $A$, $B$ and~$C$, like
\juerg{the Green type formula} $\,(A^{2\rho}v,w)=(A^\rho v,A^\rho w)$ for every $v\in\VA{2\rho}$ and $w\in\VA\rho$.

Before stating our well-posedness theorem, we prove some auxiliary results.
The first one is a separation property enjoyed by any solution under our assumptions on the data.

\Bthm
\label{Separation}
Assume \eqref{hpz} and $u\in\L2H$, and let $\soluz$ be a solution to problem \Pbl\ satisfying \Regsoluz.
Then \juerg{it holds for every $M>0$ that if $\norma\mu_\infty<M$, then}
there exists a compact interval $[\aM,\bM]\subset(a,b)$ such that
\Beq
  \aM \leq \phi \leq \bM
  \quad \aeQ.
  \label{separationM}
\Eeq
\juerg{This} interval depends only on~$f$, the initial datum~$\phiz$, and~$M$.
\Ethm

\Bdim
Notice that $\mu$ is bounded thanks to~\eqref{mubdd}.
So, we fix a constant $M$ and assume that $\norma\mu_\infty<M$.
By the assumptions \eqref{hpf} on $f$ and \eqref{hpz} on~$\phiz$,
we can choose $\aM\in(a,\az]$ and $\bM\in[\bz,b)$ such~that
\Beq
  f(z) < -M
  \quad \hbox{for \pier{all} $z\in(a,\aM)$}
  \aand
  f(z) > M
  \quad \hbox{for \pier{all} $z\in(\bM,b)$}.
  \non
\Eeq
Now, we notice that, for a.a.~$s\in(0,T)$, the value $\phi(s)$ belongs to $\VB{2\sigma}$ by~\eqref{regphi}.
Moreover, the function $z\mapsto\psi(z):=(z-\bM)^+$ is monotone and \Lip\ continuous on~$\erre$ and vanishes at the origin.
Hence, \pier{we have that} $\psi(\phi(s))\in H$ by \eqref{hpBs}.
So, we can multiply \eqref{seconda}, written at the time~$s$, by $\psi(\phi(s))$ 
and integrate \juerg{over} $(0,t)$ with respect to~$s$.
By noting that $\psi(\phiz)=0$ (since $\phiz\leq\bz\leq\bM$ \aeO), we obtain that
\Beq
  \frac 12 \, \norma{\psi(\phi(t))}^2
  + \iot \bigl( B^{2\sigma}\phi(s),\psi(\phi(s)) \bigr) \ds
  = \iot \bigl( \mu(s) - f(\phi(s)) , \psi(\phi(s)) \bigr) \ds \,.
  \non
\Eeq
Thanks to the inequality \eqref{hpBs}, 
the second term on the \lhs\ is nonnegative.
Moreover, the \rhs\ is nonpositive 
since $\psi(\phi)=0$ where $\phi\leq\bM$, and $f(\phi)\geq\mu$ whenever $\phi>\bM$.
Hence, we conclude that $\psi(\phi)=0$ \aeQ, i.e., that $\phi\leq\bM$ \aeQ.
By the same argument, with $\psi(z):=-(z+\pier{a_M})^-$, one obtains that $\phi\geq\aM$ \aeQ.
\Edim

\Bthm
\label{Contdep}
Under the assumptions \eqref{hpz} on the initial data,
problem \Pbl\ has at most one solution satisfying \Regsoluz.
Moreover, if $M>0$, $u_i\in\L2H$, $i=1,2$, and $(\mu_i,\phi_i,S_i)$ are
two corresponding solutions to \Pbl\ satisfying \Regsoluz\ and $\norma{\mu_i}_\infty<M$ for $i=1,2$,
then the estimate 
\begin{align}
  & \norma{\mu_1-\mu_2}_{\L2{\VA\rho}}
  + \norma{\phi_1-\phi_2}_{\H1H\cap\L\infty{\VB\sigma}}
  \non
  \\
  & \quad {}
  + \norma{S_1-S_2}_{\L\infty H\cap\L2{\VC\tau}}
  \leq K_M \, \norma{u_1-u_2}_{\L2H}
  \label{contdepM}
\end{align}
holds true with a constant $K_M$ that depends only on the structure of the system, the initial data, $T$ and~$M$.
\Ethm

\Bdim
We show the uniqueness at the end
and first prove the estimate \eqref{contdepM},  noting that
the assumption $\norma{\mu_i}_\infty<M$ is meaningful
since the functions $\mu_i$ are bounded due to~\eqref{mubdd}.
We apply Theorem~\ref{Separation} 
and find a compact interval $[\aM,\bM]$ contained in $(a,b)$ such~that
\Beq
  \aM \leq \phi_i \leq \bM \quad \aeQ,
  \quad \hbox{for $i=1,2$}.
  \non
\Eeq
Since $f$ is (at least) a $C^1$ function on~$(a,b)$,
it is \Lip\ continuous on~$[\aM,\bM]$.
Let $L$ be the corresponding \Lip\ constant.
After this preparation, we can start the proof.
We set, for convenience,
$u:=u_1-u_2$, $\mu:=\mu_1-\mu_2$, $\phi:=\phi_1-\phi_2$, and $S:=S_1-S_2$,
write \accorpa{prima}{terza} for both solutions and multiply the differences by
$\mu$, $\dt\phi$, and~$S$, respectively, in the inner product of~$H$.
Then, we sum up and integrate with respect to time \juerg{over}~$(0,t)$.
By noting that the terms involving $(\mu,\dt\phi)$ cancel each other,
and adding the same contributions 
$\,(1/2)\norma{\phi(t)}^2=\iot(\phi(s),\dt\phi(s))\ds\,$ 
and $\,\iot\norma{S(s)}^2\ds\,$
to both sides,
we obtain the equation
\begin{align}
  & \iot \norma{A^\rho\mu(s)}^2 \ds
  + \iot \norma{\dt\phi(s)}^2
  + \frac 12 \, \norma{\phi(t)}^2
  + \frac 12 \, \norma{B^\sigma\phi(t)}^2
  \non
  \\
  & \quad {}
  + \frac 12 \, \norma{S(t)}^2
  + \iot \norma{S(s)}^2 \ds
  + \iot \norma{C^\tau S(s)}^2 \ds
  \non
  \\
  \separa
  & = \iot \Bigr( P(\phi_1(s))(S_1(s)-\mu_1(s)) - P(\phi_2(s))(S_2(s)-\mu_2(s)) , \mu(s)-S(s) \bigr) \ds
  \non
  \\
  & \quad {}
  - \iot \bigl( f(\phi_1(s)) - f(\phi_2(s)) , \dt\phi(s) \bigr) \ds 
  \non
  \\
  & \quad {}
  + \iot \bigl( \phi(s) , \dt\phi(s) \bigr) \ds 
  + \iot \bigl( u(s) , \juerg{S(s)} \bigr) \ds 
  + \iot \norma{S(s)}^2 \ds \,.
  \non
\end{align}
If we term $I$ the first integral on the \rhs, 
apply the Young inequality to the next three terms,
use the \Lip\ continuity of~$f$, and rearrange,
we deduce that
\begin{align}
  & \iot \norma{A^\rho\mu(s)}^2 \ds
  + \juerg{\frac 12} \iot \norma{\dt\phi(s)}^2
  + \frac 12 \, \norma{\phi(t)}_\Bs^2
  + \frac 12 \, \norma{S(t)}^2
  + \iot \norma{S(s)}_\Ct^2 \ds
  \non
  \\
  & \leq I
  + (L^2 + 1) \iot \norma{\phi(s)}^2 \ds 
  + \juerg{\frac 12}\iot \norma{u(s)}^2 
  + \juerg{\frac 32}\iot \norma{S(s)}^2 \ds \,.
  \non
\end{align}
Now, we rewrite $I$ as the sum of two terms.
The first of these is nonpositive, since $P$ is nonnegative, 
and we estimate the other one recalling that $P'$ is bounded, because $P$ is \Lip\ continuous.
By applying the \Holder\ inequality, and recalling~\eqref{embconst}, we have for every $\delta>0$ that
\begin{align}
  I\,& = \iot \bigl( P(\phi_1(s))(S(s)-\mu(s)) , \mu(s) - S(s) \bigr) \ds
  \non
  \\
  &  \quad + \iot \bigl( (P(\phi_1(s)))-P(\phi_2(s))) (S_2(s)-\mu_2(s)) , \mu(s)-S(s) \bigr) \ds
  \non
  \\
  %\separa
  & \leq\, \sup|P'| \iot \norma{\phi(s)}_4 \, \norma{S_2(s)-\mu_2(s)}_4 \, \juerg{\|S(s)-\mu(s)\|} \ds
  \non
  \\
  & \leq\, \delta \iot \bigl( \norma{S(s)} + \juerg{\norma{\mu(s)}} \bigr)^2 \ds
  \non
  \\
  &\quad {}
  +\, \frac {\sup|P'|^2 \, C_*^4}{4\delta}
    \iot \bigl( \norma{S_2(s)}_\Ct + \norma{\mu_2(s)}_\Ar \bigr)^2 \norma{\phi(s)}_\Bs^2 \ds \,.
  \non
\end{align}
We notice that the function $s\mapsto\bigl( \norma{S_2(s)}_\Ct + \norma{\mu_2(s)}_\Ar \bigr)^2$
belongs to \pier{$L^\infty(0,T)$}, thanks to the regularity \eqref{regmu} and \eqref{regS} of $\mu_2$ and~$S_2$, respectively.
Therefore, by choosing $\,\delta>0\,$ small enough, and applying the Gronwall lemma,
we obtain the estimate \eqref{contdepM} with a constant $\,K_M\,$ 
whose \juerg{dependence on the data} agrees with that specified in the statement.

We now come back to uniqueness.
As both $u_i\in\L2H$ and the corresponding solutions $(\mu_i,\phi_i,S_i)$ are arbitrary in the above argument 
(since no restriction on $M$ is made), \juerg{we conclude that \,$(\mu_1,\phi_1,S_1)=(\mu_2,\phi_2,S_2)\,$
if $u_1=u_2$, which shows the uniqueness for the solution to problem \Pbl.
With this}, the proof is complete.
\Edim

Finally, we can state our well-posedness and stability result.
Here, and later on in this paper, we use the  notation
\Beq
  \BR := \graffe{u\in\L2H :\ \norma u_{\L2H} < R},
  \label{defBR}
\Eeq
where $R$ is a positive real parameter.

\Bthm
\label{Wellposedness}
Under the assumptions \eqref{hpz} on the initial data $\phiz$ and~$\Sz$,  
\juerg{problem \Pbl\ has for every $u\in\L2H$ a unique solution $\soluz$ that satisfies} \Regsoluz.
In particular, $\mu$ is bounded.
Moreover, for every $R>0$, there exist a constant $K_1(R)$ 
and a compact interval $[\aR,\bR]\subset(a,b)$,
which depend only on the structure of the system, the initial data, $T$ and $R$, such that
both the estimate
\begin{align}
  & \norma\mu_{\L\infty{\VA{2\rho}}}
  + \norma\mu_\infty+ \norma{P(\phi)(S-\mu)}_{\L2H}
  \non
  \\
  & \quad {}
  + \norma\phi_{\W{1,\infty}H \cap \H1{\VB\sigma} \cap \L2{\VB{2\sigma}}}
  + \norma{f(\phi)}_{\L\infty H} 
  \non
  \\
  & \quad {}
  + \norma S_{\H1H \cap \L\infty{\VC\tau} \cap \L2{\VC{2\tau}}}
    \non
  \\
  & \leq K_1(R) 
  \label{stimasoluz}
\end{align}
and the separation property
\Beq
  \aR \leq \phi \leq \bR
  \quad \aeQ
  \label{separation}
\Eeq
hold true for every $u\in\BR$ and the corresponding solution $\soluz$.
Finally, if $R>0$, $u_i\in\BR$, $i=1,2$, and $(\mu_i,\phi_i,S_i)$ are the corresponding solutions,
\juerg{then} the estimate 
\begin{align}
  & \norma{\mu_1-\mu_2}_{\L2{\VA\rho}}
  + \norma{\phi_1-\phi_2}_{\H1H\cap\L\infty{\VB\sigma}}
  \non
  \\
  & \quad {}
  + \norma{S_1-S_2}_{\L\infty H\cap\L2{\VC\tau}}
  \leq K_2(R) \, \norma{u_1-u_2}_{\L2H}
  \label{contdep}
\end{align}
holds true with a constant $K_2(R)$ that depends only on the structure of the system, the initial data, $T$, and~$R$.
\Ethm

\Bdim
Uniqueness follows from Theorem~\ref{Contdep}.
Let us come to the existence of a solution and to the estimates \eqref{stimasoluz} and~\eqref{contdep}.
However, we do not give a complete proof.
Indeed, one \juerg{can} adapt the arguments of~\cite{CGS23}
on account of Remark~\ref{HpCGS23},
and we briefly explain the reason \juerg{for this}.
The procedure used there is based on the Yosida regularization of the nonlinearity~$f$, 
a time discretization of the regularized system,
and the derivation of suitable a priori estimates,
and the same line \juerg{of argumentation can be followed in our situation.}
Here is the main remark: in~\cite{CGS23}, some estimates for~$\mu$ have been derived from estimates of~$\dt\mu$,
and this term is missing in \eqref{prima}, in contrast to~\eqref{Iprima}.
In the present case, an estimate of the norm of~$\mu$, e.g., in $\L2H$,
can be deduced from an estimate of $A^\rho\mu$ in the same space 
as shown in Remark~\ref{Firsteigenvalue},
by using the assumption \eqref{hpeigen} on the first eigenvalue $\lambda_1$ of~$A$.
Hence, we do not repeat the arguments of \cite{CGS23} with \juerg{the corresponding} modifications.
However, for the reader's convenience,
we sketch the formal proofs of the estimates
that would be obtained step by step in the \juerg{rigorous} procedure 
in order to prove the existence of a solution.
We assume $u\in\BR$ from the very beginning,
so that these estimates eventually lead to~\eqref{stimasoluz} as well.

In order to simplify notation, we use the same symbol $c$ without any subscript
for possibly different constants that depend only 
on the structure of our system, the initial data and~$T$, but neither on~$u$ nor on $R$\,.
Moreover, the~symbol $c_R$ stands for (possibly different) constants
that depend on the constant~$R$, in addition, but still not on~$u$.
So, it is understood that the actual values of such constants may vary 
from line to line and even in the same chain of inequalities.
Notice that the notations used for constants we want to refer to
(like, e.g., those used in~\eqref{hpbelow})
are different.

\step
First a priori estimate

We test \eqref{primav}, \eqref{secondav} and \eqref{terzav}, written at the time $\,s$, 
by~$\mu(s)$, \juerg{$\dt\phi(s)$}, and~$S(s)$, respectively, in the scalar product of~$H$.
Then we sum up and integrate over~$(0,t)$, where $t\in(0,T)$ is arbitrary,
noting that the terms involving the product $\mu\,\dt\phi$ cancel each other.
By also adding $\,|\Omega|C_2\,$ to both sides (see \eqref{hpbelow}), we obtain the identity
\begin{align}
  & \iot \norma{A^\rho\mu(s)}^2 \ds
  + \iot \norma{\dt\phi(s)}^2 \ds
  + \frac 12 \, \norma{B^\sigma\phi(t)}^2 
  + \iO (F(\phi(t)) + C_2)
  \non
  \\
  \separa
  & \quad {}
  + \frac 12 \, \norma{S(t)}^2
  + \iot \norma{C^\tau S(s)}^2 \ds
  + \intQt P(\phi) (S-\mu)^2
  \non
  \\
  \separa
  & {} = \frac 12 \, \norma{B^\sigma\phiz}^2 
  + \iO (F(\phiz) + C_2)
  + \frac 12 \, \norma\Sz^2 
  + \iot \bigl( u(s) , S(s) \bigr) \ds \,.
  \non
\end{align}
Recalling the consequence \eqref{forequiv} of~\eqref{hpeigen},
\pier{taking \eqref{hpbelow} into account}
%eliminating~$C_2$ in the norms by means of the triangle inequality 
and applying the Gronwall lemma, 
we conclude that
\begin{align}
  & \norma\mu_{\L2{\VA\rho}}
  + \norma\phi_{\H1H\cap\L\infty{\VB\sigma}}
  \non
  \\[1mm]
  & \quad {}
  + \norma S_{\L\infty H\cap\L2{\VC\tau}}
  + \norma{P^{1/2}(\phi)(S-\mu)}_{\L2H}
  \leq c_R \,. 
  \qquad
  \label{primastima}
\end{align}

\step
\juerg{Consequence}

\juerg{We can also \pier{infer} that
\Beq
  \norma{P(\phi)(S-\mu)}_{\L2H}^2
  + \norma{F(\phi)}_{\L\infty\Luno}
  \leq c_R \,.
  \non
\Eeq
Indeed, the estimate is right for the first term, since $P$ is bounded by~\eqref{hpP}; for the bound 
of the second term,
we can argue as in \cite[Section 4.4]{CGS23}.
} 
%On the contrary, that for the latter is not correct at the present moment.
%Indeed, the true estimate should be performed on the discrete scheme of the regularized system
%(here one uses the uniform boundedness from below with respect to $\lambda$ mentioned in Remark~\ref{HpCGS23})
%and it does not seems to be conserved in the limit as $\lambda$ tends to zero
%(because of weak types of convergence in spaces constructed on~$Q$).
%Indeed, just the \lhs\ of the inequality
%\Beq
%  T^{-1} \, \norma{F(\phi)}_{\LQ1}
%  \leq \norma{F(\phi)}_{\L\infty\Luno}
%  \non  
%\Eeq
%will be bounded.

\step 
Second a priori estimate

We would like to test \eqref{primav} by $\dt\mu$ even though $\dt\mu$ does not appear in the equation
(in contrast to~\eqref{Iprima}).
In fact, the estimate we derive here by a formal procedure
should be performed rigorously at the level of the discete scheme,
which can contain that time derivative multiplied by a viscosity coefficient 
that tends to zero at some point of the procedure.
So, we test \eqref{prima} and \eqref{terza} \juerg{formally} by $\dt\mu$ and~$\dt S$, respectively.
At the same time, we \juerg{formally} differentiate \eqref{seconda} with respect to time 
and test the resulting equality by~$\dt\phi$.
Then, we sum up and integrate with respect to time\pier{,} as usual.
Since the terms involving the product $\dt\mu\,\dt\phi$ cancel each other,
we obtain the identity
\begin{align}
  & \frac 12 \, \norma{A^\rho\mu(t)}^2
  + \frac 12 \, \norma{\dt\phi(t)}^2
  + \iot \norma{B^\sigma\dt\phi(s)}^2 \ds
  \non
  \\
  & \quad {}
  + \iot \norma{\dt S(s)}^2 \ds
  + \frac 12 \, \norma{C^\tau S(t)}^2
  \non
  \\
  & {} = \frac 12 \, \norma{\juerg{A^\rho\mu}(0)}^2
  + \frac 12 \, \norma{\dt\phi(0)}^2
  + \frac 12 \, \norma{C^\tau\Sz}^2
  \non
  \\
  & \quad {}
  + \iot \bigl( u(s) , \juerg{\dt S(s)} \bigr) \ds 
  - \iot \bigl( f'(s) \dt\phi(s) , \dt\phi(s) \bigr) \ds 
  \non
  \\
  & \quad {}
  + \iot \bigl( P(\phi(s))(S(s)-\mu(s)) , \dt\mu(s) - \dt S(s) \bigr) \ds \,.
  \label{persecondastima}
\end{align}
The integral containing $u$ \juerg{can be handled using Young's inequality}, 
and the one involving $f'$ is easily treated 
\pier{using} the second inequality \pier{in} \eqref{hpbelow} and~\eqref{primastima}.
We \pier{postpone the estimate of the last integral}
and  first deal with the initial values appearing on the \rhs\ of \eqref{persecondastima}. 
By using the initial conditions for $\phi$ and~$S$,
we write \pier{\eqref{prima} and \eqref{seconda}} at the time $t=0$
in the following way:
\Beq
  \bigl( A^{2\rho} + P(\phiz) \bigr) \mu(0)
  = P(\phiz) \Sz
  \aand
  \dt\phi(0) 
  = \mu(0) - B^{2\sigma} \phiz - f(\phiz) \,.
  \label{perstimaz}
\Eeq
Since $P$ is nonnegative, by multiplying the first \juerg{identity in}  \eqref{perstimaz} by~$\mu(0)$,
we have that (see~\eqref{forequiv})
\Beq
  \bigl( P(\phiz)\Sz , \mu(0) \bigr)
  = \bigl( ( A^{2\rho} + P(\phiz) ) \mu(0) , \mu(0) \bigr)
  \geq \norma{A^\rho \mu(0)}^2
  \geq \lambda_1^{2\rho} \, \norma{\mu(0)}^2\,,
  \non
\Eeq
whence
\Beq
  \norma{\mu(0)}
  \leq \lambda_1^{-2\rho} \, \norma{P(\phiz)\Sz}\,
  \juerg{\leq \,c}, \quad 
   \pier{ \norma{A^\rho \mu(0)}^2
  \leq \norma{\mu(0)} \, \norma{P(\phiz)\Sz}\,
  \leq \,c}
  \label{stimamuz}
\Eeq
and, on account of~\eqref{hpz}, we also deduce from the second identity in \eqref{perstimaz} that
\Beq
  \norma{\dt\phi(0)}
  \leq c + \norma{B^{2\sigma}\phiz} + \norma{f(\phiz)}
  \leq c \,.
  \non
\Eeq
Finally, we deal with the last integral on the \rhs\ of \eqref{persecondastima}, 
which we term $I$ for brevity.
We perform an integration by parts in time, recall that $P'$ is bounded by~\eqref{hpP},  and invoke \eqref{stimamuz}.
By recalling that $\norma\cpto_p$ denotes the norm in $\Lx p$, we then have
\begin{align}
  I\,& =\, - \iot \!\!\! \iO P(\phi) (S-\mu) \dt(S-\mu)
  \non
  \\
  & 
  = \,- \frac 12 \iO P(\phi(t)) \bigl( S(t) - \mu(t) \bigr)^2
  + \frac 12 \iO P(\phiz) \bigl( \Sz - \mu(0) \bigr)^2
  \non
  \\
  & \quad {}
  + \frac 12 \iot \!\!\! \iO P'(\phi) \dt\phi (S-\mu)^2
  \non
  \\
  & \leq\, c 
  + c \iot \norma{\dt \phi(s)}_2 \bigl(
    \norma{S(s)}_4^2 
    + \norma{\mu(s)}_4^2
  \bigr) \ds \,.
  \non 
\end{align}
Finally, we notice that
\Beq
  \ioT \bigl( \norma{S(s)}_4^2 + \norma{\mu(s)}_4^2 \bigr) \ds
  \leq c_R\,,
  \non
\Eeq
by \eqref{primastima} and some of the embeddings in \eqref{embeddings}. 
\juerg{This allows us to apply Gronwall's lemma, whence we conclude that}
\Beq
  \norma\mu_{\L\infty{\VA\rho}}
  + \norma{\dt\phi}_{\L\infty H\cap\L2{\VB\sigma}}
  + \norma S_{\H1H\cap\L\infty{\VC\tau}}
  \leq c_R \,.
  \label{secondastima}
\Eeq

\step
Third a priori estimate

By \pier{taking $v=\mu(t)$ in \eqref{primav}}  and recalling that $P$ is nonnegative and bounded, 
we obtain the following inequality \aat
\Beq
  \norma{A^\rho\mu(t)}^2 
  \leq \bigl( P(\phi(t))S(t) - \dt\phi(t) , \mu(t) \bigr)
  \leq c \, \bigl( \norma S_{\L\infty H} + \norma{\dt\phi}_{\L\infty H} \bigr) \norma{\mu(t)} \,.
  \non
\Eeq
By accounting for Remark~\ref{Firsteigenvalue}, we deduce that
\Beq
  \norma\mu_{\L\infty{\VA\rho}} \leq c_R \,.
  \label{terzastima}
\Eeq
In particular, the norm of $\mu$ in $\L\infty H$ is bounded by \pier{some constant $c_R$.}
Since the same holds for $S$ due to \eqref{primastima},
we infer that
\Beq
  \norma{P(\phi)(S-\mu)}_{\L\infty H} \leq c_R \,.
  \label{stimaP}
\Eeq

\step
Consequence

By comparison in \eqref{prima},
we deduce that
\Beq
  \norma{A^{2\rho}\mu}_{\L\infty H}
  \leq \norma{P(\phi)(S-\mu)}_{\L\infty H}
  + \norma{\dt\phi}_{\L\infty H} .
  \non
\Eeq
Combining this with \eqref{stimaP} and \eqref{secondastima}, we conclude that
\Beq
  \norma\mu_{\L\infty{\VA{2\rho}}}
  \leq c_R \,.
  \label{stimamumax}
\Eeq
Then, the first embedding in \eqref{embeddings} yields that $\mu$ is bounded
(as \juerg{claimed} in the statement) and that
\Beq
  \norma\mu_\infty 
  \leq c_R \,.
  \label{stimamuinfty}
\Eeq
By comparison in \eqref{terza}, we also deduce that
\Beq
  \norma S_{\L2{\VC{2\tau}}}
  \leq c_R \,.
\label{stimaSmax} 
\Eeq
No better estimate for $S$ is available, since $u\in\L2H$, only.

\step\juerg{Fourth} a priori estimate

We recall that Remark~\ref{HpCGS23} provides a splitting of $f$ as $f_1+f_2$, 
with $f_1$ monotone and vanishing at the origin and $f_2$ \Lip\ continuous.
So, we can write \eqref{seconda} \aat\ in the form
\Beq
  B^{2\sigma}\phi(t) + f_1(\phi(t))
  = \mu(t) - \dt\phi(t) - f_2(\phi(t)),
  \non
\Eeq
and test this identity by $f_1(\phi(t))$.
More precisely, in the correct argument $f_1$ is replaced by its Yosida regularization,
which is monotone and \Lip\ continuous and vanishes at the origin,
and the equation itself is replaced by a scheme, which is obtained by 
\pier{discretizing} time differentiation
and for which the analogue of $\phi(t)$ belongs to $\VB{2\sigma}$.
Hence, assumption \eqref{hpBs} can actually be applied.
Here, we formally apply it to the above identity with $v=\phi(t)$ and $\psi=f_1$.
We obtain that
\begin{align}
  & \norma{f_1(\phi(t))}^2
  \,\leq\, \bigl( B^{2\sigma}\phi(t) + f_1(\phi(t)) , f_1(\phi(t)) \bigr)
  \,=\, \bigl( \mu(t) - \dt\phi(t) - f_2(\phi(t)) , f_1(\phi(t)) \bigr)
  \non
  \\[0.5mm]
  & \leq \,\norma{\mu-\dt\phi-f_2(\phi)}_{\L\infty H} \, \norma{f_1(\phi(t))} \,.
  \non
\end{align}
On account of the previous estimates, and by a comparison in \eqref{seconda},
we conclude that
\Beq
  \norma{f_1(\phi)}_{\L\infty H} + \norma\phi_{\L\infty{\VB{2\sigma}}}
  \leq c_R \,.
  \label{stimafphi} 
\Eeq

\step 
Conclusion

This concludes the formal proof of the existence part of Theorem~\ref{Wellposedness}
and of estimate~\eqref{stimasoluz}.
As already said, in the \juerg{rigorous} argument
the above bounds are established for the solution to an approximating problem,
and one has to perform some limiting procedure.
The estimates provide convergence of weak and weak-star type.
However, even strong convergence in $\L2H$ for the approximations of $\phi$ and $S$ is obtained.
Indeed, the embeddings $\VB\sigma\subset H$ and $\VC\tau\subset H$ are compact due to~\eqref{hpABC},
so that one can apply the Aubin--Lions lemma
(see, e.g., \cite[Thm.~5.1, p.~58]{Lions}).
Therefore, the nonlinear terms can be correctly managed.

\step Separation

Let us come to estimate~\eqref{separation}.
This is a trivial consequence of the above estimate and Theorem~\ref{Separation}.
Indeed, this theorem can be applied with $M:=c_R+1$, where $c_R$ is the constant that appears in~\eqref{stimamuinfty}.
The corresponding compact interval $[\aR,\bR]$ of the statement 
is nothing but the interval $[\aM,\bM]$ \pier{considered} in \eqref{separationM},
which depends only on the structure of the system, the initial data, $T$, and~$R$.

\step
Continuous dependence

Also \eqref{contdep} is a trivial consequence of a fact already proved, namely, of Theorem~\ref{Contdep}.
Indeed, if $u_i\in\BR$, $i=1,2$, \juerg{then}
the $L^\infty$ bound for the corresponding $\mu_i$ is ensured by \eqref{stimasoluz},
and Theorem~\ref{Contdep} can be applied with \gianni{\hbox{$M=K_1(R)+1$}}.
Hence, we can take as $K_2(R)$ the constant $K_M$ that appears in~\eqref{contdepM}.
\juerg{Also this} constant depends only on the structure of the system, the initial data, $T$ and~$R$.
\Edim

\juerg{%
\Brem
The existence part of Theorem 2.6 is closely connected to the existence result proved in
\cite[Theorem~3.4]{CGS24}, where, however, no statement
concerning separation or uniqueness was proved. For purposes of control theory, however, it is
indispensable to have uniqueness, since otherwise no control-to-state operator can be defined,
and this seems to be  available only under the assumptions made here.
\Erem
}%

%%%%%%%%%%%%%%%%%%%%%%%%%%%%%%%%%%%%%%%%%%%%%%%%%%%%%%%%%%%%%%%%%%%%%%%%

\section{The control-to-state mapping}
\label{FRECHET}
\setcounter{equation}{0}

The results of the previous section ensure that we can correctly define a control-to-state mapping
to be used in the control problem \juerg{under investigation}.
Taking into account that the cost functional \juerg{to be minimized} 
depends only on the components $\phi$ and $S$ of the solution corresponding to a given~$u$,
we~set
\begin{align}
  & \calY_1 := \L2{\VA\rho} , \quad
  \calY_2 :=\H1H\cap\L\infty{\VB\sigma}, 
  \non
  \\[0.5mm]
  &\calY_3 := \pier{\C0 H{}}\cap\L2{\VC\tau},
  \aand
  \calY := \calY_2 \times \calY_3,
  \label{defY}
\end{align}
and define
\Beq
  \calS_i : \L2H \to \calY_i \,, \quad i=1,2,3, \
  \aand
  \calS : \L2H \to \calY ,
  \non
\Eeq
by setting, for $u\in\L2H$,
\begin{align}
  & \hbox{$\calS_1(u):=\mu$, \ $\calS_2(u):=\phi$, \ $\calS_3(u):=S$, \ and \ $\calS(u):=(\phi,S)$},
  \non
  \\[0.5mm]
  & \hbox{where $\soluz$ is the solution to \Pbl\ corresponding to $u$}.
  \qquad
  \label{defcalS}
\end{align}
More precisely, we need to consider the restriction of these maps to $\BR$ for any given radius~$R>0$.
The choice of the space $\,\calY\,$ mainly is due to the following fact:
the inequality \eqref{contdep} implies~that
\Beq
  \norma{\calS(u_1)-\calS(u_2)}_\calY \leq K_2(R) \, \norma{u_1-u_2}_{\pier{\L2H}}
  \quad \hbox{for every $u_1\,,u_2\in\BR$}\,.
  \label{lipcontostate}
\Eeq
A very important consequence of the separation property \eqref{separation} and of the regularity of $f$ 
\pier{ensured by} \eqref{hpF}  is the following global boundedness condition:
\Beq
  \norma{f^{(k)}(\calS_2(u))}_\infty
  \leq K_3(R)
  \quad \hbox{for $k=0,1,2$, and every $u\in\BR$},
  \label{gb}
\Eeq
where $\,\juerg{K_3(R)}\,$ depends only on the structure of the system, the initial data, $T$, and~$R$.

The \Frechet\ differentiability of the maps $\calS$ 
is \pier{strictly} related to the properties of the \juerg{linearized} problem we introduce \juerg{now. To this end,
we} fix $\ub\in\L2H$.
The linearized system associated with $\ub$ and the variation $h\in\L2H$ is the following:
\begin{align}
  & \dt\xi + A^{2\rho}\eta
  = P(\phib) (\zeta-\eta)
  + P'(\phib) \, \xi \, (\Sb-\mub)\,,
  \label{primal}
  \\[0.5mm]
  & \dt\xi + B^{2\sigma} \xi + f'(\phib) \, \xi 
  = \eta\,,
  \label{secondal}
  \\[0.5mm]
  & \dt\zeta + C^{2\tau} \zeta
  = - P(\phib) (\zeta-\eta)
  - P'(\phib) \, \xi \, (\Sb-\mub)
  + h\,,
  \label{terzal}
  \\[0.5mm]
  & \xi(0) = 0 \aand \zeta(0) = 0\,, 
  \label{cauchyl}
\end{align}
\Accorpa\Pbll primal cauchyl
where $\mub:=\calS_1(\ub)$, $\phib:=\calS_2(\ub)$, and $\Sb:=\calS_3(\ub)$.
We also write the weak formulation of \accorpa{primal}{terzal}:
the identities
\begin{align}
  & \bigl( \dt\xi(t) , v \bigr)
  + \bigl( A^\rho \eta(t) , A^\rho v \bigr) 
  \non
  \\
  & = \bigl( P(\phib(t)) (\zeta(t)-\eta(t)) , v \bigr)
  + \bigl( P'(\phib(t)) \, \xi(t) \, (\Sb(t)-\mub(t)) , v \bigr)\,,
  \label{primalv}
  \\[2mm]
  & \bigl( \dt\xi(t) , v \bigr)
  + \bigl( B^\sigma\xi(t) , B^\sigma v \bigr)
  + \bigl( f'(\phib(t)) \, \xi(t) , v \bigr)
  \non
  \\
  & = \bigl( \eta(t) , v \bigr)\,,
  \label{secondalv}
  \\[2mm]
  & \bigl( \dt\zeta(t) , v \bigr)
  + \bigl( C^\tau \zeta(t) , C^\tau v \bigr)
  \non
  \\
  & = - \bigl( P(\phib(t)) (\zeta(t)-\eta(t)) , v \bigr)
  \non
  \\
  & \quad {}
  - \bigl( P'(\phib(t)) \, \xi(t) \, (\Sb(t)-\mub(t)) , v \bigr)
  + \bigl( h(t) , v \bigr)\,,
  \label{terzalv}
\end{align}
\juerg{have to hold true} for every $v\in\VA\rho$, $v\in\VB\sigma$, and $v\in\VC\tau$, respectively, and \aat.
\juerg{We have the following results.}

\Bthm
\label{Wellposednessl}
\juerg{Suppose that the assumptions \eqref{hpz} on the initial data of problem \Pbl\ are fulfilled, and}
let $\ub\in\L2H$ and $h\in\L2H$.
Then the linearized problem \Pbll\ has a unique solution $\soluzl$ satisfying the regularity requirements
\begin{align}
  & \eta \in \L2{\VA{2\rho}}\,,
  \label{regeta}
  \\[0.5mm]
  & \xi \in \H1H \cap \L\infty{\VB\sigma} \cap \L2{\VB{2\sigma}}\,,
  \label{regxi}
  \\[0.5mm]
  & \zeta \in \H1H \cap \L\infty{\VC\tau} \cap \L2{\VC{2\tau}}.
  \label{regzeta}
\end{align}
\Accorpa\Regsoluzl regeta regzeta
Moreover, if $R>0$ and $\ub\in\BR$, \juerg{then} this solution satisfies the estimate
\Beq
  \norma{(\xi,\zeta)}_\calY \leq K_4(R) \, \norma h_{\L2H}\,,
  \label{stimal}
\Eeq
where the constant $K_4(R)$ depends only on the structure of the system \accorpa{prima}{terza}, 
the initial data $\phiz$ and~$\Sz$, $T$, and~$R$.
\Ethm

\Bdim
We notice that the coefficients $P(\phib)$, $P'(\phib)$, $f'(\phib)$, \juerg{as well as~$\mub$,} 
are bounded functions.
Moreover, if $\ub$ belongs to some~$\BR$, \juerg{then}
the $L^\infty$~bounds are uniform, i.e., they \juerg{just} depend on $R$ and not on~$\ub$.
On the contrary, $\Sb$~might be unbounded. 
However, as stated in Theorem~\ref{Wellposedness}, it is smooth.
So, the linear system is not worse than the nonlinear one
and can be solved by the same argument (which we do not repeat here)
based on time discretization that has been used in~\cite{CGS23} 
(see also~\cite{CGS19} for the linearized system associated with the Cahn--Hilliard equations).
This first leads to a solution to the variational problem \accorpa{primalv}{terzalv}
and then to the strong formulation~\accorpa{primal}{terzal}.
However, we perform at least some formal estimates 
that can justify both the regularity asserted in the statement and the validity of estimate~\eqref{stimal}.
To this end, we fix $R>0$ and assume that $\ub\in\BR$ at once.

Also in this section, i.e., in this proof and later on, 
we adopt a convention on the constants similar to the one used in the previous section:
$c$~stands for possibly different constants depending only on the structure, the data, and $T$,
while the notation $c_R$ \juerg{indicates an additional dependence} on~$R$.

\step
First a priori estimate

We formally test \accorpa{primalv}{terzalv} by $\eta$, $\dt\xi$, and~$\zeta$, respectively,
sum up and integrate over~$(0,t)$.
Moreover, we add the same quantities $\,(1/2)\norma{\xi(t)}^2=\iot(\xi(s),\dt\xi(s))\ds\,$ and $\,\iot\norma{\zeta(s)}^2\,$
to both sides of the resulting identity, 
in order to recover the full norms in $\VB\sigma$ and $\VC\tau$ on the \lhs.
Also in this case a cancellation occurs, and we have that
\begin{align}
  & \iot \norma{A^\rho\eta(s)}^2 \ds
  + \iot \norma{\dt\xi(s)}^2 \ds
  + \frac 12 \, \norma{\xi(t)}_\Bs^2
  + \frac 12 \, \norma{\zeta(t)}^2
  + \iot \norma{\zeta(s)}_\Ct^2 \ds
  \non
  \\
  \separa
  & = \iot \bigl( P(\phib(s)) (\zeta(s)-\eta(s)) , \juerg{\eta(s)} - \zeta(s) \bigr) \ds
  \non
  \\
  & \quad {}
  + \iot \bigl( P'(\phib(s)) \, \xi(s) \, (\Sb(s)-\mub(s)) , \juerg{\eta(s)} - \zeta(s) \bigr) \ds
  \non
  \\
  & + \iot \bigl( \xi(s) - f'(\phib(s)) \, \xi(s) , \dt\xi(s) \bigr) \ds
  + \iot \bigl( h(s) + \zeta(s) , \zeta(s) \bigr) \ds \,.
  \label{perprimastimal}
\end{align}
The first integral on the \rhs\ is nonpositive, \juerg{while} the next one, which we term~$I$, needs some treatment.
By using the \Holder\ inequality, two of the inequalities \eqref{embconst},
Remark~\ref{Firsteigenvalue}, and the Young inequality,
we obtain that
\begin{align}
  I\,& \leq \,c \iot \norma{\xi(s)}_2 \, \norma{\Sb(s)-\mub(s)}_4 \, \norma{\eta(s)-\zeta(s)}_4 \ds
  \non
  \\
  & \leq\, c \iot \norma{\xi(s)} \, \bigl( \norma{\Sb(s)}_\Ct + \norma{\mub(s)}_\Ar \bigr)
    \, \bigl( \norma{\eta(s)}_\Ar + \norma{\zeta(s)}_\Ct \bigr) \ds
  \non
  \\
  & \leq\, \frac 12 \iot \norma{A^\rho\eta(s)}^2 \ds
  + \frac 12 \iot \norma{\zeta(s)}_\Ct^2 \ds
  + c \iot \bigl( \norma{\Sb(s)}_\Ct^2 + \norma{\mub(s)}_\Ar^2 \bigr) \norma{\xi(s)}^2 \ds\,,
  \non
\end{align}
where we notice that the function $\,s\mapsto\norma{\Sb(s)}_\Ct^2 + \norma{\mub(s)}_\Ar^2\,$
belongs to $L^\infty(0,T)$ and that its norm is bounded by a constant like in~\eqref{stimasoluz},
due to Theorem~\ref{Wellposedness} applied to~$\ub$.
By treating the last terms of \eqref{perprimastimal} using the Schwarz and Young inequalities,
and applying Gronwall's lemma, we conclude~that
\begin{align}
  & \norma\eta_{\L2{\VA\rho}}
  + \norma\xi_{\pier{\H1 H\cap\L\infty {\VB\sigma}}}
  + \norma\zeta_{\L\infty H\cap\L2{\VC\tau}}
  \non
  \\[0.5mm]
  & \leq c_R \, \norma h_{\L2H} \,.
  \label{primastimal}
\end{align}

\step
Second a priori estimate

We estimate the \rhs\ of \eqref{primal}.
On account of the embeddings~\eqref{embeddings}, we have \aet\ that
\begin{align}
  & \norma{P(\phib) (\zeta-\eta) + P'(\phib) \, \xi \, (\Sb-\mub)}
  \,\leq\, c \, \bigl( \norma\zeta + \norma\eta \bigr)
  + c \,\norma\xi_4 \, \norma{\Sb-\mub}_4 
  \non
  \\
  & \leq \,c \, \bigl( \norma\zeta + \norma\eta \bigr)
  + c \, \norma\xi_\Bs \,\bigl( \norma\Sb_\Ct + \norma\mub_\Ar \bigr).
  \non
\end{align}
On account of \eqref{primastimal}, we conclude that
\Beq
  \norma{P(\phib) (\zeta-\eta) + P'(\phib) \, \xi \, (\Sb-\mub)}_{\L2H} 
  \leq c_R \, \norma h_{\L2H} \,.
  \label{stimaRHS}
\Eeq
Since \eqref{primastimal} also yields an estimate for $\dt\xi$,
a~(formal) comparison in \eqref{primal} allows us to conclude that
\Beq
  \norma{A^{2\rho}\eta}_{\L2H}
  \leq c_R \, \norma h_{\L2H} \,, 
  \quad \hbox{i.e.,} \quad
  \norma\eta_{\L2{\VA{2\rho}}}
  \leq c_R \, \norma h_{\L2H} \,.
  \label{secondastimal}
\Eeq

\step
Third a priori estimate

We test \eqref{terzal} by $\dt\zeta$ and integrate in time, as usual.
On account of~\eqref{stimaRHS} and the Young inequality, we obtain that
\Beq
  \iot \norma{\dt\zeta(s)}^2 \ds
  + \frac 12 \, \norma{C^\tau \zeta(t)}^2
  \leq \frac 12 \, \iot \norma{\dt\zeta(s)}^2 \ds
  + c_R \, \norma h_{\L2H}^2 \,.
  \non
\Eeq
We thus deduce that
\Beq
  \norma{\dt\zeta}_{\L2H} + \norma\zeta_{\L\infty{\VC\tau}}
  \leq c_R \, \norma h_{\L2H} \,.
\Eeq
Now that $\dt\zeta$ is estimated, 
a~comparison in \pier{\eqref{terzal}} provides a bound for $C^{2\tau}\zeta$.
Hence, we conclude that
\Beq 
  \norma\zeta_{\H1H\cap\L\infty{\VC\tau}\cap\L2{\VC{2\tau}}}
  \leq c_R \, \norma h_{\L2H} \,. 
  \label{terzastimal}
\Eeq

This ends the list of the formal estimates 
and formally leads to a strong solution satisfying~\eqref{stimal}.
Even though uniqueness formally follows by taking $h=0$,
we remark that it can be proved rigorously.
Indeed, by assuming the regularity \Regsoluzl, 
the procedure used to obtain the above estimates is justified.
\Edim

\Bthm
\label{Frechet}
Assume \eqref{hpz} for the initial data of problem \Pbl.
Then the control-to-state mapping $\,\calS\,$ defined in \eqref{defcalS}
is \Frechet\ differentiable at every point in $\L2H$.
More precisely, if $\ub\in\L2H$ and $h\in\L2H$, \juerg{then}
the value $(D\calS)(\ub)[h]$ of the \Frechet\ derivative $(D\calS)(\ub)$ in the direction $h$
is given by \pier{the pair $(\xi, \zeta)$, where $\soluzl$ is the solution} 
to the linearized problem \Pbll\ associated with $\ub$ and~$h$.
\Ethm

\Bdim
Fix any $\ub\in\L2H$, and let $\soluzb$ be the corresponding state.
For every $h\in\L2H$, let $(\muh,\phih,\Sh)$ be the state corresponding to $\ub+h$.
Finally, let $\soluzl$ be the solution to the linearized problem \Pbll\ associated with $\ub$ and~$h$.
We set, for convenience,
\Beq
  \etah := \muh - \mub - \eta \,, \quad
  \xih := \phih - \phib - \xi
  \aand
  \zetah := \Sh - \Sb - \zeta \,.
  \label{defdiff}
\Eeq
According to the definitions of differentiability and derivative in the sense of \Frechet,
we have to prove that the (linear) map $h\mapsto(\xi,\zeta)$ is continuous from $\L2H$ into~$\calY$
and that there exist a real number $\overline h>0$ and
a function $\Lambda:(0,\overline h)\to\erre$ satisfying
\Beq
  \norma{(\xih,\zetah)}_\calY
  \leq \Lambda(\norma h_{\L2H})
  \aand
  \lim_{s\searrow 0} \frac {\Lambda(s)} s = 0 \,.
  \label{frechet}
\Eeq
The first fact is ensured by \eqref{stimal} once $R$ is chosen larger than $\norma\ub_{\L\infty H}$.
Hence, we fix $R>\norma\ub_{\L2H}$ once and for all.
As for the construction of~$\Lambda$, we set $\overline h:=R-\norma\ub_{\L2H}$,
and we assume that $\norma h_{\L2H}<\overline h$.
This implies that $\ub$ and $\ub+h$ belong to~$\BR$,
so that Theorem~\ref{Wellposedness} can be applied to both of them.
We \juerg{thus} derive uniform estimates for the corresponding states,
\juerg{hence} for the coefficients of the corresponding linearized systems.
This \juerg{entails} uniform estimates for the corresponding solutions.
In order to establish~\eqref{frechet},
we observe that $(\etah,\xih,\zetah)$ satisfies the regularity properties
\begin{align}
  & \etah \in \L\infty{\VA{2\rho}}\,,
  \non
  \\[0.5mm]
  & \xih \in \pier{\H1H \cap \L\infty{\VB\sigma} \cap \L2{\VB{2\sigma}}}\,,
  \non
  \\[0.5mm]
  & \zetah \in \H1H \cap \L\infty{\VC\tau} \cap \L2{\VC{2\tau}}\,,
  \non
\end{align}
and solves the problem
(in~a strong form, i.e., the equations are satisfied \aeQ, 
since all the contributions are $L^2$ functions)
\begin{align}
  & \dt\xih + A^{2\rho} \etah 
  = \Qh1\,,
  \label{primah}
  \\[0.5mm]
  & \dt\xih + B^{2\sigma}\xih + \Qh2
  = \etah\,,
  \label{secondah}
  \\[0.5mm]
  & \dt\zetah + C^{2\tau} \zetah
  = - \Qh1\,,
  \label{terzah}
  \\[0.5mm]
  & \etah(0) = 0\, , \quad
  \xih(0) = 0\,,
  \aand
  \zetah(0)=0\,,
  \label{cauchyh}
\end{align}
where $\Qh1$ and $\Qh2$ are defined by
\begin{align}
  & \Qh1 := P(\phih) (\Sh-\muh) 
  - P(\phib) (\Sb-\mub)
  - P(\phib) (\zeta-\eta)
  - P'(\phib) \, \xi \, (\Sb-\mub) \,,
  \non
  \\[0.5mm]
  & \Qh2 := f(\phih) - f(\phib) - f'(\phib) \, \pier{\xi}\,.
  \non
\end{align}
It is convenient to rewrite the functions $\Qh i$
by accounting for the Taylor expansions of $P$ and~$f$.
\pier{Usig the formula with integral remainder, it}
is immediately checked~that
\begin{align}
\label{repQ2}
\Qh2 = f'(\phib) \xih + \Rh2 \, (\phih-\phib)^2,
  \quad \hbox{where} \quad
  \Rh2 := \int_0^1 (1-\theta) f''(\phib+\theta(\phih-\phib)) \, d\theta\,,
\end{align}
while it is more complicated to find a convenient representation of~$\Qh1$.
However, \juerg{simple algebraic manipulations
show} that
\begin{align}
  & \Qh1 = P(\phib) (\zetah-\etah)
  + ( P(\phih) - P(\phib)) [(\Sh-\Sb) - (\muh-\mub)]
  \non
  \\
  & \qquad {}
  + P'(\phib) (\Sb-\mub) \, \xih
  + (\Sb-\mub) \, \Rh1 \, (\phih-\phib)^2\,,
  \label{repQ1}  
\end{align}
where
$$ \Rh1 := \int_0^1 (1-\theta) P''(\phib+\theta(\phih-\phib)) \, d\theta \,.
$$
Notice that, by~\eqref{gb}, both $\Rh1$ and $\Rh2$ are bounded uniformly with respect to~$h$:
\Beq
  \norma{\Rh1}_\infty + \norma{\Rh2}_\infty
  \leq c_R \,.
  \label{stimeRh}
\Eeq

After this preparation, we start estimating. \juerg{To this end, we} test \eqref{primah}, \eqref{secondah},
 and \eqref{terzah}, by $\etah$, $\dt\xih$, and~$\zetah$, respectively.
Then, we sum up and integrate in time.
There is a usual cancellation.
By adding the same contributions to both sides similarly as in the previous proof, we obtain
\begin{align}
  & \iot \norma{A^\rho\etah(s)}^2 \ds
  + \iot \norma{\dt\xih(s)}^2 \ds
  + \frac 12 \, \norma{\xih(t)}_\Bs^2
  \non
  \\
  & \quad {}
  + \frac 12 \, \norma{\zetah(t)}^2
  + \iot \norma{\zetah(s)}_\Ct^2 \ds
  \non
  \\
  & = \iot \bigl( \Qh1(s) , \etah(s)-\zetah(s) \bigr) \ds
  - \iot \bigl( \Qh2(s) , \dt\xih(s) \bigr) \ds 
  \non
  \\
  & \quad {}
  + \iot \bigl( \xih(s) , \dt\xih(s) \bigr) \ds 
  + \iot \norma{\zetah(s)}^2 \ds \,.
  \label{perfrechet}
\end{align}
We have to \pier{estimate} only the integrals involving $\Qh1$ and $\Qh2$.
The first term produces four integrals, termed $I_j$ for $j=1,\dots,4$ for brevity,
which correspond to the four summands, in that order, of~\eqref{repQ1}.
Clearly, $I_1$~is nonpositive.
As for~$I_2$, we use the \Holder\ inequality and the embeddings~\eqref{embeddings}
as well as the estimate~\pier{\eqref{contdep}} applied with $u_1=\ub+h$ and $u_2=\ub$.
Hence, by omitting the integration variable~$s$ to shorten the lines, we have for every~$\delta>0$
\begin{align}
   I_2\,& \leq \,c \iot \norma{\phih-\phib}_4
    \, \bigl( \juerg{\norma{\Sh-\Sb} + \norma{\muh-\mub}} \bigr)
    \, \bigl( \norma{\etah}_4 + \norma{\zetah}_4 \bigr) \ds
  \non
  \\
  & \leq\, c \iot \norma{\phih-\phib}_\Bs
    \, \bigl( \juerg{\norma{\Sh-\Sb} + \norma{\muh-\mub}} \bigr)
    \, \bigl( \norma{\etah}_\Ar + \norma{\zetah}_\Ct \bigr) \ds
  \non
  \\
  & \leq\, \delta \iot  \bigl( \norma{\etah}_\Ar^2 + \norma{\zetah}_\Ct^2 \bigr) \ds
  \non
  \\
  & \quad {}
  + \frac c\delta \, \norma{\phih-\phib}_{\L\infty{\VB\sigma}}^2 \bigl( \norma{\Sh-\Sb}_{\L2H}^2 + \norma{\muh-\mub}_{\L2H}^2 \bigr)
  \non
  \\
  & \leq\, \delta \iot  \bigl( \norma{\etah}_\Ar^2 + \norma{\zetah)}_\Ct^2 \bigr) \ds
  + \frac {c_R} \delta \, \norma h_{\L2H}^4 \,.
  \non
\end{align}
Next, by also accounting for \eqref{stimasoluz} applied to $\mub$ and~$\Sb$, we similarly have that
\begin{align}
   I_3\,& \leq\, c \iot \norma{\Sb-\mub}_4 \, \norma\xih_4\, \bigl( \juerg{\norma\etah + \norma\zetah} \bigr) \ds
  \non
  \\
  & \leq\, c \iot \bigl( \norma\Sb_\Ct + \norma\mub_\Ar \bigr) \norma\xih_\Bs \bigl( \norma\etah + \norma\zetah \bigr) \ds
  \non
  \\
  & \leq\, \delta \iot \norma\etah^2 \ds
  + \juerg{{c_R}\bigl(1+\delta^{-1}\bigr)} \, \iot \bigl( \juerg{\|\xi^h\|^2_{B,\sigma}+ \norma\zetah^2} \bigr) \ds 
  \non
\end{align}
and for the fourth contribution, \pier{thanks} to~\eqref{stimeRh}, we obtain that
\begin{align}
   I_4\,& \leq\, c_R \iot \norma{\Sb-\mub}_4 \, \norma{\phih-\phib}_4^2 \, \norma{\etah-\zetah}_4 \ds
  \non
  \\
  & \leq\, c_R \iot \bigl( \norma\Sb_\Ct + \norma\mub_\Ar \bigr) \, \norma{\phih-\phib}_\Bs^2 \, \bigl( \norma\etah_\Ar + \norma\zetah_\Ct \bigr) \ds
  \non
  \\
  & \leq\, \delta \iot \bigl( \norma\etah_\Ar^2 + \norma\zetah_\Ct^2 \bigr) \ds
  + \frac {c_R}\delta \, \norma h_{\L2H}^4 \,.
  \non
\end{align}
Finally, we estimate the term involving~$\Qh2$ by accounting for \eqref{repQ2} and \eqref{stimeRh} in this~way:
\begin{align}
  & - \iot ( \Qh2 , \dt\xih ) \ds
  \,=\, - \iot \bigl( \pier{f'(\phib)} \xih + \Rh2 \, (\phih-\phib)^2 , \dt\xih \bigr) \ds  
  \non
  \\
  &\le\, c_R \iot \norma \xih \, \norma{\dt\xih} \ds
  + c_R \iot \norma{\phih-\phib}_4^2 \, \juerg{\norma{\dt\xih}} \ds
  \non
  \\
  & \leq\, \delta \iot \norma{\dt\xih}^2 \ds
  + \frac {c_R}\delta \, \iot \norma\xih^2 \ds
  + \frac {c_R}\delta \, \norma h_{\L2H}^4 \,.
  \non
\end{align}
By treating the last two terms of \eqref{perfrechet} in a trivial way,
recalling all the inequalities derived above, choosing $\,\delta>0\,$ small enough, 
and applying the Gronwall lemma,
we conclude~that
\Beq
  \juerg{\|\etah\|_{L^2(0,T;V_A^\rho)}}
  + \norma\xih_{\H1H\cap\L\infty{\VB\sigma}}
  + \norma\zetah_{\L\infty H\cap\L2{\VC\tau}}
  \leq c_R \, \norma h_{\L2H}^2 \,.
  \non
\Eeq
If we term $C_R$ the value of the constant $c_R$ of the last inequality,
then we obtain \eqref{frechet} with $\Lambda$ defined on $(0,\overline h)$ by $\Lambda(s):=C_R\, s^2$.
This completes the proof.
\Edim

%%%%%%%%%%%%%%%%%%%%%%%%%%%%%%%%%%%%%%%%%%%%%%%%%%%%%%%%%%%%%%%%%%%%%%%%

\section{The control problem}
\label{CONTROL}
\setcounter{equation}{0}

As \juerg{announced} in the Introduction, the main aim of this paper 
is the discussion of a control problem for the state system studied in the previous sections.
For this problem, we assume that
\begin{align}
  & \kappa_i \geq 0, \quad \hbox{for $i=1,\dots,5$}, \quad
  \phi_Q \,,\, S_Q \in \LQ2,
  \aand
  \phi_\Omega \,,\, S_\Omega \in \Ldue \,,
  \qquad
  \label{hpcost}
  \\
  & \umin \,,\, \umax \in \LQ\infty,
  \aand
  \umin \leq \umax \quad \aeQ.
  \label{hpumm}
\end{align}
\Accorpa\HPcontrol hpcost hpumm
Then, the cost functional $\,\calJ\,$ and the set $\Uad$ of the admissible controls are defined by
\begin{align}
  & \calJ(u,\phi,S)
  := \juerg{\frac{\kappa_1}2} \intQ |\phi-\phi_Q|^2
  + \juerg{\frac{\kappa_2}2} \iO |\phi(T) - \phi_\Omega|^2
  \non
  \\
  & \quad {}
  + \juerg{\frac{\kappa_3}2} \intQ |S-S_Q|^2
  + \juerg{\frac{\kappa_4}2} \iO |S(T) - S_\Omega|^2
  + \juerg{\frac{\kappa_5}2} \intQ |u|^2\,,
  \label{cost}
  \\
  & \Uad := \graffe{ u\in\L2H : \ \umin \leq u \leq \umax \ \aeQ}\,,
  \label{defUad}
\end{align}
and the control problem is the following:
\begin{align}
  & \hbox{Minimize} \quad \calJ(u,\phi,S) \quad \hbox{under the constraints that} \quad
  u \in \Uad
  \quad \hbox{and}
  \non
  \\
  & \mbox{$\soluz$ is the solution to \Pbl\ corresponding to $u$}.
  \qquad
  \label{controlpbl}
\end{align}
For the above problem, we prove the existence of an optimal control,
and we derive the first-order necessary conditions for optimality.
This involves an adjoint problem for which we prove a well-posedness result.
We recall that the control-to-state mapping $\calS$ is defined in~\eqref{defcalS}
and state our first result.

\Bthm
\label{Optimum}
Under the assumptions \HPcontrol, the control problem has at least one solution, that is,
there \juerg{is some}\, $\ub\in\Uad\,$ satisfying the following condition:
for every $\juerg{v}\in\Uad$ we have that $\,\calJ(\ub,\phib,\Sb)\leq\calJ(v,\phi,S)$,
where $(\phib,\Sb)=\calS(\ub)$ and $(\phi,S)=\calS(v)$.
\Ethm

\Bdim
Since $\Uad$ is nonempty, the infimum of $\calJ$ under the constraints given in \eqref{controlpbl}
is a well-defined real number $d\geq0$,
and we can pick a minimizing sequence~$\graffe{\un}\subset\Uad$.
\juerg{Hence, denoting by $(\mun,\phin,\Sn)$ the state corresponding to~$\un$ for $n\in\enne$, we have that
$\calJ(\un,\phin,\Sn)\to d$ as $n\to\infty$.}
Since $\Uad$ is bounded and closed in~$\L2H$ (in~fact, it is even bounded and closed in~$\LQ\infty$),
we can assume that
\Beq
  \un \to \ub 
  \quad \hbox{weakly in $\L2H$}
  \label{convun}
\Eeq
for some $\ub\in\Uad$.
Moreover, we can choose some $R>0$ such that $\Uad\subset\BR$.
Therefore, we can apply Theorem~\ref{Wellposedness} to $\un$
and deduce that $(\mun,\phin,\Sn)$ satisfies the estimate~\eqref{stimasoluz}, as well as
the separation and global boundedness properties \eqref{separation} and~\eqref{gb}, for all $n\in\enne$.
Hence, for a subsequence indexed again by $\,n$, we have~that
\begin{align}
  & \mun \to \mub
  \quad \hbox{weakly star in $\L\infty{\VA{2\rho}}$}\,,
  \label{convmun}
  \\[0.5mm]
  & \phin \to \phib
  \quad \hbox{weakly star in $\W{1,\infty}H\cap\H1{\VB\sigma}\cap{\L2{\VB{2\sigma}}}$}\,,
  \qquad
  \label{convphin}
  \\[0.5mm]
  & \Sn \to \Sb
  \quad \hbox{weakly star in $\H1H\cap\L\infty{\VC\tau}\cap\L2{\VC{2\tau}}$}\,.
  \label{convSn}
\end{align}
It follows that the initial conditions \eqref{cauchy} are satisfied by the limiting pair~$(\phib,\Sb)$.
Moreover, thanks to the compact embedding $\VB\sigma\subset H$ ensured by~\eqref{hpABC},
\pier{and consequently of $\H1{\VB\sigma}$ into $\L2H$, we}
%we can apply the Aubin--Lions lemma
%(see, e.g., \cite[Thm.~5.1, p.~58]{Lions}) and 
deduce that
\Beq
  \phin \to \phib
  \quad \hbox{strongly in $\L2H$}.
  \non
\Eeq
Since $f$ and $P$ are \Lip\ continuous in $[\aR,\bR]$,
we \pier{also infer that}
\Beq
  f(\phin) \to f(\phib)
  \aand
  P(\phin) \to P(\phib),
  \quad \hbox{strongly in $\L2H$}.
  \non
\Eeq
It follows that $(\mub,\phib)$ solves~\eqref{seconda}.
From the above strong convergence and the weak convergence of $\{\mun\}$ and $\{\Sn\}$ at least in $\L2H$,
we deduce that $\{P(\phin)(\Sn-\mun)\}$ converges to $P(\phib)(\Sb-\mub)$ weakly in~$\LQ1$.
Hence, the limiting triplet $\soluzb$ satisfies equations \eqref{prima} and \eqref{terza} as well,
i.e., $\soluzb$ is the state corresponding to the control~$\ub$.
On the other hand, we have that 
\Beq
  \calJ(\ub,\phib,\Sb)
  \,\leq\, \liminf_{n\nearrow\infty} \calJ(\un,\phin,\Sn)
  = d\,,
  \non
\Eeq
by semicontinuity.
We conclude that $\ub$ is an optimal control.
\Edim

The rest of the section is devoted to the derivation of the first-order necessary conditions for optimality.
Hence, we fix an optimal control $\ub\in\Uad$ and the corresponding $\soluzb$ once and for all. \pier{If we introduce}
\juerg{the so-called reduced cost functional $\tilde\calJ$ by setting}
\Beq
  \tilde\calJ(u) := \calJ(u,\calS_2(u),\calS_3(u))
  \quad \hbox{for $u\in\L2H$}\,,
  \non
\Eeq
we immediately find \juerg{from the convexity of $\Uad$} that
the \Frechet\ derivative $(D\tilde\calJ)(\ub)\in\calL(\L2H;\erre)$ must satisfy
\Beq
  (D\tilde\calJ)(\ub)[v-\ub] \geq 0
  \quad \hbox{for every $v\in\Uad$},
  \non
\Eeq
provided that it exists.
But this is the case 
due to the obvious differentiability of the quadratic functional $\calJ$
and the differentiability of~the operator~$\calS$,
which takes its values in~$\calY\subset(\C0H)^2$.
Hence, by accounting for the full statement of Theorem~\ref{Frechet},
we can even apply the chain rule and rewrite the above inequality~as
\begin{align}
  & \kappa_1 \intQ (\phib - \phi_Q) \xi
  + \kappa_2 \iO (\phib(T) - \phi_\Omega) \xi(T)
    + \kappa_3 \intQ (\Sb - S_Q) \zeta
 \non\\   
  &+ \kappa_4 \iO (\Sb(T) - S_\Omega) \zeta(T)
    + \kappa_5 \intQ \ub (v-\ub)
  \,\,\juerg{\geq}\,\, 0
  \quad \hbox{for every $v\in\Uad$},
  \label{badnc}
\end{align}
where $\xi$ and $\zeta$ are the components of the solution $\soluzl$
to the linearized system \Pbll\ associated with $\ub$ and $h=v-\ub$.

As usual in control problems, a condition of this sort is not satisfactory,
since it requires to solve the linearized problem for infinitely many choices of $h\in\L2H$,
because $v$ is arbitrary in~$\Uad$.
Therefore, we have to eliminate $\xi$ and $\zeta$ from~\eqref{badnc},
which can be done by introducing and solving a proper adjoint problem.
This is a backward-in-time problem \juerg{for the adjoint state variables $(q,p,r)$} that formally reads
as follows:
\begin{align}
  & A^{2\rho} q - p + P(\phib) (q-r) 
  \,=\, 0\,, 
  \label{primaa}
  \\[0.5mm]
  & - \dt(q+p) + B^{2\sigma} p  
  + f'(\phib) \, p - P'(\phib) (\Sb-\mub) (q-r)
  \,=\, \kappa_1 (\phib-\phi_Q)\,,
  \qquad
  \label{secondaa}
  \\[0.5mm]
  & - \dt r + C^{2\tau} r - P(\phib) (q-r) 
  \,=\, \kappa_3(\Sb-S_Q)\,, 
  \label{terzaa}
  \\[0.5mm]
  & (q+p)(T) \,=\, \kappa_2 (\phib(T)-\phi_\Omega)
  \aand 
  r(T) \,=\, \kappa_4 (\Sb(T) - S_\Omega) \,.
  \label{cauchya}
\end{align}
\Accorpa\Pbla primaa cauchya
However, in order to give this system a proper meaning according to the regularity that we will prove,
we need some preliminaries.
First, due to the density of $\VB\sigma$ in~$H$,
we can identify $\,H\,$ \juerg{with} a subspace of the dual space $\,\VB{-\sigma}:=(\VB\sigma)^*\,$ of~$\,\VB\sigma\,$
in \juerg{such a way}  that $\<v,w>=(v,w)$ for every $v\in H$ and $w\in\VB\sigma$,
where $\<\cpto,\cpto>$ denotes the duality pairing between $\VB{-\sigma}$ and~$\VB\sigma$.
Now, thanks to the obvious formula $(B^{2\sigma} v,w)=(B^\sigma v,B^\sigma w)$,
which holds for every $v\in\VB{2\sigma}$ and $w\in\VB\sigma$
and, \juerg{owing to the above identification, can also be read in the form $\<B^{2\sigma} v,w>=(B^\sigma v,B^\sigma w)$,}
one can extend the operator $B^{2\sigma}:\VB{2\sigma}\to H$
to a continuous linear operator, still termed $B^{2\sigma}$, from $\VB\sigma$ to $\VB{-\sigma}$
by means of the above formula, namely,
\Beq
  \< B^{2\sigma} v , w > = ( B^\sigma v , B^\sigma w )
  \quad \hbox{for every $v,w\in\VB\sigma$}.
  \label{ext}
\Eeq
At this point, it is meaningful to postulate the following regularity for the adjoint variables:
\begin{align}
  & q \in  \L\infty{\VA{2\rho}}\,,
  \label{regq}
  \\[0.5mm]
  & p \in \L2{\VB\sigma}
  \aand
  \dt(q+p) \in \L2{\VB{-\sigma}}\,,
  \label{regp}
  \\[0.5mm]
  & r \in \H1H \cap \L\infty{\VC\tau} \cap \L2{\VC{2\tau}}\,.
  \label{regr}
\end{align}
\Accorpa\Regsoluza regq regr
Indeed, then all of the equations, as well as the final conditions, have a precise meaning,
by also accounting for the properties of the other ingredients which we recall for the reader's convenience:
$P(\phib)$, $P'(\phib)$, $f'(\phib)$, \pier{and $\mub$ are bounded}, and
$\Sb\in\H1H\cap\L\infty{\VB\sigma}$, whence, in particular, $\Sb\in\L\infty{\Lx4}$.
However, we also consider a variational formulation of the adjoint system,
which makes sense in a much weaker regularity setting for~$\soluza$, namely,
\Beq
  q \in  \L\infty{\VA\rho}, \quad
  p \in \L2{\VB\sigma},
  \aand
  r \in \H1H \cap \L2{\VC\tau}.
  \label{regav}
\Eeq
We require that
\begin{align}
  & \ioT \bigl\{ \bigl( A^\rho q , A^\rho v \bigr)
  - ( p,v ) + \bigl( P(\phib) (q-r) ,v ) \bigr\} \ds
  = 0 
  \non
  \\
  & \quad \hbox{for every $v\in\L2{\VA\rho}$},
  \label{primaav}
  \\
  & \ioT \bigl\{
    (q+p , \dt v) 
    + ( B^\sigma p , B^\sigma v ) + \bigl( f'(\phib) \, p - P'(\phib) (\Sb-\mub) (q-r) , v \bigr)
  \bigr\} \ds
  \non
  \\
  & {}
  = \ioT (g_1,v) \ds + \bigl( g_2 , v(T) \bigr)
  \non
  \\
  & \quad \hbox{for every $v\in\H1H\cap\L2{\VB\sigma}$ vanishing at $t=0$},
  \label{secondaav}
  \\
  & \ioT \bigl\{
    ( -\dt r , v)
    + ( C^\tau r , C^\tau v) - \bigl( P(\phib) (q-r) , v \bigr)
  \bigr\} \ds
  \non
  \\
  & {}
  = \ioT (g_3 , v ) \ds 
  \quad \hbox{for every $v\in\L2{\VC\tau}$},
  \label{terzaav}
  \\
  & r(T) = g_4 ,
  \label{cauchyav}
\end{align}
\Accorpa\Pblav primaav cauchyav
where we have \juerg{introduced the abbreviating notation}
\Beq
  g_1 := \kappa_1 (\phib-\phi_Q) , \
  g_2 := \kappa_2 (\phib(T)-\phi_\Omega) , \
  g_3 := \kappa_3(\Sb-S_Q) , \
  g_4 := \kappa_4 (\Sb(T) - S_\Omega) .
  \quad
  \label{abbrev}
\Eeq
Also for brevity, and in order to shorten the exposition,
we have omitted the integration time variable termed~$s$.
We will do the same in the following.

Clearly, \Regsoluza\ and \Pbla\ imply \eqref{regav} and \Pblav.
In fact, these problems are equivalent.
The proof given below makes use of the Leibniz rule
proved in \cite[Lem.~4.5]{CGS18}
(and well known under slightly different assumptions),
which we here state as a lemma.

\Blem
\label{Leibniz}
Let $(\calV,\calH,\calV^*)$ be a Hilbert triplet,
and assume that
\Beq
  y \in \H1\calH \cap \L2\calV
  \aand
  z \in \H1{\calV^*} \cap \L2\calH \,.
  \label{hpleibniz}
\Eeq
Then the function $t\mapsto(y(t),z(t))_{\calH}$
is absolutely continuous on~$[0,T]$, 
and its derivative is given~by
\Beq
  \frac d{dt} \, (y,z)_{\calH}
  = (y',z)_{\calH}
  + {}_{\calV^*}\< z',y>_{\calV}
  \quad \aet,
  \label{leibniz}
\Eeq
\juerg{where $(\,\cdot\,,\,\cdot\,)_{\cal H}$ and ${}_{\calV^*}\langle\,\cdot\,,\,\cdot\,\rangle_{\calV}$ denote the
inner product in $\calH$ and the dual pairing between $\calV^*$ and $\calV$, respectively.}
\Elem

\Blem
\label{Equivalence}
Assume \juerg{that \eqref{regav} and \Pblav\ are valid}.
Then \Regsoluza\ and \Pbla\  hold true as well.
\Elem

\Bdim
We first notice that \eqref{primaav} implies the pointwise variational inequality
\Beq
  \bigl( A^\rho q , A^\rho v \bigr)
  = \bigl( p - P(\phib) (q-r) , v \bigr)
  \quad \hbox{for every $v\in\VA\rho$ and \aet}.
  \non
\Eeq
On the other hand, the conditions
$w\in\VA\rho$, $g\in H$, and $(A^\rho w,A^\rho v)=(g,v)$ for every $v\in\VA\rho$,
imply that $w\in\VA{2\rho}$ and $A^{2\rho}w=g$,
as one immediately sees by using the spectral representation.
Hence, we obtain \eqref{regq} and~\eqref{primaa}.
The same argument can be used to deduce that
$r$ belongs to $\L2{\VC{2\tau}}$ and solves~\eqref{terzaa},
since even $\dt r$ belongs to $\L2H$ by assumption.
The last condition $r\in\L\infty{\VC\tau}$ in \eqref{regr} then follows from interpolation.

Much more work has to be done for the second equations.
First, for the same test functions $v$ as in~\eqref{secondaav}, we deduce~that
\Beq
  \ioT \bigl\{
    (q+p , \dt v) 
    + \< B^{2\sigma} p , v >
  \bigr\} \ds
  = \ioT (g,v) \ds + \bigl( g_2 , v(T) \bigr)\,,
  \label{forsec}
\Eeq
where, for brevity, we have set 
\Beq
  g := g_1 - f'(\phib) \, p + P'(\phib) (\Sb-\mub) (q-r) \,.
  \non
\Eeq
We immediately infer that
\begin{align}
  & \pier{{}\left|\ioT \bigl( q+p , \dt v \bigr) \ds \right|{}}
  \leq \norma p_{\L2{\VB\sigma}} \, \norma v_{\L2{\VB\sigma}}
  \non
  \\
  & \quad {}
  + \norma g_{\L2H} \, \norma v_{\L2H}
  + \norma{g_2} \, \norma{v(T)} \,.
  \non
\end{align}
In particular, we have for some constant $c>0$ that
\Beq
  \pier{{}\left|\ioT \bigl( q+p , \dt v \bigr) \ds \right|{}}
  \leq c \, \norma v_{\L2{\VB\sigma}}
  \quad \hbox{for every $v\in C^\infty_c(0,T;\VB\sigma)$}.
  \non  
\Eeq
This exactly means that $\dt(q+p)\in(\L2{\VB\sigma})^*=\L2{\VB{-\sigma}}$.
Thus, we can replace the expression $\,(q+p,\dt v)\,$ by $\,-\<\dt(q+p),v>\,$ in~\eqref{forsec}, provided that $v\in C^\infty_c(0,T;\VB\sigma)$.
The variational equation we obtain \juerg{is just \eqref{secondaa} understood} in the sense of~$\VB{-\sigma}$.

It remains to derive the first of the final conditions~\eqref{cauchya}.
To this end, we also assume that $v(0)=0$ and \pier{exploit \eqref{forsec}} once more.
Moreover, we can apply Lemma~\ref{Leibniz}
with $\calV=\VB\sigma$, $\calH=H$, $y=v$, and $z=q+p$,
since $v\in\H1H\cap\L2{\VB\sigma}$ and $q+p\in\H1{\VB{-\sigma}}\cap\L2H$.
Finally, we account for the already proved equation \eqref{secondaa}.
We then obtain that
\begin{align}
  & \ioT (g,v) \ds + \bigl( g_2,v(T) \bigr)
  = \ioT (q+p,\dt v) \ds + \ioT (B^\sigma p,B^\sigma v) \ds
  \non
  \\
  & = \ioT \bigl\{ - \< \dt(q+p) , v > + \< B^{2\sigma} p , v > \bigr\} \ds
  + \< (q+p)(T) , v(T) >
  \non
  \\
  & = \ioT (g,v) \ds 
  + \< (q+p)(T) , v(T) > \,.
  \non
\end{align}
Therefore, we have that $\bigl((q+p)(T),v(T)\bigr)=\bigl(g_2,v(T)\bigr)$ for every $v$ with the required properties,
and the desired final condition obviously follows.
\Edim

So, we can choose between the strong form \Pbla\ and the weak formulation~\Pblav,
according to our convenience, in proving a well-posedness result, which is our next goal.
We prepare the existence part by introducing a Faedo--Galerkin scheme with viscosity
that looks like an approximation of \Pbla.
We recall \eqref{eigen} on the eigenvalues and the eigenvectors of the operators and set,
for every integer~$n>1$, 
\Beq
  \VAn := \mathop{\rm span} \graffe{e_1,\dots,e_n} , \quad
  \VBn := \mathop{\rm span} \graffe{e'_1,\dots,e'_n},
  \aand
  \VCn := \mathop{\rm span} \graffe{e''_1,\dots,e''_n} \,.
  \non
\Eeq
Then, we look for a triplet $\soluzn$ satisfying
\Beq
  \qn \in \H1\VAn , \quad
  \pn \in \H1\VBn,
  \aand
  \rn\in\H1\VCn,
  \label{regFG}
\Eeq
and solving the system
\begin{align}
  & \bigl( - \textstyle\frac 1n \, \dt\qn + A^{2\rho} \qn - \pn + P(\phib) (\qn-\rn) , v \bigr)
  \non
  \\
  & \quad {} 
  = 0
  \quad \hbox{for every $v\in\VAn$ and \aet},
  \label{priman}
  \\
  \separa
  & \bigl( - \dt(\qn+\pn) + \bigl( B^{2\sigma} \pn + f'(\phib) \pn - P'(\phib) (\Sb-\mub) (\qn-\rn) , v \bigr)
  \non
  \\
  & \quad {} 
  = ( g_1 , v )
  \quad \hbox{for every $v\in\VBn$ and \aet},
  \label{secondan}
  \\
  \separa
  & \bigl( - \dt\rn + C^{2\tau} \rn - P(\phib) (\qn-\rn) , v \bigr)
  \non
  \\
  & \quad {} 
  = ( g_3 , v )
  \quad \hbox{for every $v\in\VCn$ and \aet},
  \label{terzan}
\end{align}
as well as the final conditions
\begin{align}
  & \bigl( \qn(T) , v \bigr) = 0 , \quad
  \bigl( (\qn+\pn)(T) , v \bigr)  = \bigl( g_2 , v \bigr),
  \aand
  \bigl( \rn(T) , v \bigr) = \bigl( g_4 , v \bigr),
  \non
  \\
  & \quad \hbox{for every $v\in\VAn$, $v\in\VBn$, and $v\in\VCn$, respectively}.
  \label{cauchyn}
\end{align}
\Accorpa\Pbln priman cauchyn

The following result holds true.

\Bprop
\label{WellposednessFG}
The system \accorpa{priman}{cauchyn} has a unique solution $\soluzn$ satisfying the conditions \eqref{regFG}.
\Eprop

\Bdim
The requirements \eqref{regFG} mean that
\Beq
  \qn(t) = \somma j1n \qn_j(t) \, \juerg{e_j} \,, \quad
  \pn(t) = \somma j1n \pn_j(t) \, \juerg{e'_j}\,,
  \aand
  \rn(t) = \somma j1n \rn_j(t) \, \juerg{e''_j}\,,
  \non
\Eeq
\aat\ and some functions $\qn_j \,,\, \pn_j \,,\, \rn_j \in H^1(0,T)$.
Moreover, an equivalent system is obtained by taking  for $\,i=1,\dots,n\,$ just 
$\,v=e_i$\,, $v=e'_i$\,, and $v=e''_i$, in the three variational equations, respectively.
Hence, \accorpa{priman}{terzan} becomes an ODE system having the column vectors
$\,q_n:=(\qn_j)$, $p_n:=(\pn_j)$, and $\,r_n:=(\rn_j),$ as unknowns.
This system reads as follows:
\begin{align}
  & \somma j1n \Bigl\{
    - \frac 1n (e_j,e_i) \frac d{dt} \qn_j 
    + \lambda_j^{2\rho} (e_j,e_i) \qn_j
    - (e'_j,e_i) \pn_j
  \non
  \\[-2mm]
  & \quad {}
    + \bigl( P(\phib) e_j,e_i \bigr) \qn_j
    - \bigl( P(\phib) e''_j,e_i \bigr) \rn_j
  \Bigr\} 
  = 0\,,
  \non
  \\
  & \somma j1n \Bigl\{
    - (e_j,e'_i) \frac d{dt} \, \qn_j
    - (e'_j,e'_i) \frac d{dt} \, \pn_j
    + (\lambda_j')^{2\sigma} (e'_j,e'_i) \pn_j
    + \bigl( f'(\phib) e'_j,e'_i \bigr) \pn_j
  \non
  \\[-2mm]
  & \quad {}
    - \bigl( P'(\phib)(\Sb-\mub) e_j,e'_i \bigr) \qn_j
    + \bigl( P'(\phib)(\Sb-\mub) e''_j,e'_i \bigr) \rn_j
  \Bigr\}
  = ( g_1 , e'_i )\,,
  \non
  \\
  & \somma j1n \Bigl\{
    - (e''_j,e''_i) \frac d{dt} \, \rn_j
    + (\lambda''_j)^{2\tau} (e''_j,e''_i) \rn_j
  \non
  \\[-2mm]
  & \quad {}
    - \bigl( P(\phib) e_j,e''_i \bigr) \qn_j
    + \bigl( P(\phib) e''_j,e''_i \bigr) \rn_j
  \Bigr\}
  = ( g_3 , e''_i )\,,
  \non
\end{align}
where the index $i$ runs over $\graffe{1,\dots,n}$ in all of the equations,
which are understood to hold \aet.
Thus, thanks to the orthogonality conditions in~\eqref{eigen}, it takes the form
\begin{align}
  & - \frac 1n \, q_n' 
  + M_1 \, q_n + M_2 \, p_n + M_3 \, r_n
  = 0\,,
  \non
  \\[0.5mm]
  & M_4 \, q_n' - p_n'
  + M_5 \, q_n + M_6 \, p_n + M_7 \, r_n 
  = b'_n\,,
  \non
  \\[0.5mm]
  & - r_n' 
  + M_8 \, q_n + M_9 \, r_n
  = b''_n\,,
  \non
\end{align}
for some (possibly time dependent, but bounded) $(n\times n)$ matrices~$\,M_k$, $k=1,\dots,9$, 
and column vectors $\,b'_n\,,\,b''_n\in\L2{\erre^n}$.
Therefore, one can solve the first equation for $\,q_n'\,$ and replace $\,q_n'\,$ in the second one
by the resulting expression.
At the same time, one multiplies the first equation by $n$ and keeps the third one as it is.
This procedure leads to an equivalent system of the form $\,-y'+My=b\,$ 
for some matrix $\,M\in\L\infty{\erre^{3n\times 3n}}\,$ and some vector $\,b\in\L2{\erre^{3n}}\,$
in the unknown $y\in\H1{\erre^{3n}}$ obtained by rearranging the triplet $\,(q_n\,,p_n\,,r_n)\,$ as a $3n$-column vector.
On the other hand, the final conditions \eqref{cauchyn} provide a final condition for~$y$.
Hence, \juerg{standard results for ODEs show the unique solvability.}
\Edim

At this point, we are ready to solve the adjoint problem. \juerg{We need, however, the following additional compatibility
condition:}

\begin{equation}
\label{compat}
\juerg{\mbox{It holds }\,\,\,\kappa_4 \,S_\Omega \in V_C^\tau}.
\end{equation} 

\Brem
\juerg{The compatibility condition \eqref{compat} is satisfied if either $\kappa_4=0$ or $S_\Omega\in V_C^\tau$. 
Obviously,
$\kappa_4=0$ means that we do not have a tracking of the solution variable $S$ at the final time $T$; while this is not
desirable, it is not too much of a restriction, since one is rather interested in monitoring the final tumor fraction
$\phi(T)$ than $S(T)$. 
On the other hand, the assumption $S_\Omega\in V_C^\tau$ is not overly restrictive in view of 
the fact that $S\in H^1(0,T;H)\cap L^2(0,T;V_C^{2\tau})$, whence it follows that $S\in C^0([0,T];V_C^\tau)$ by
continuous embedding, and thus $S(T)\in V_C^\tau$; assuming the same regularity for $S_\Omega$ is certainly not unreasonable.}  
\Erem

\juerg{We have the following result.}
\Bthm
\label{Wellposednessa}
\juerg{Suppose that also \eqref{compat} is fulfilled. Then the} adjoint system \Pbla\ has a unique solution satisfying the regularity conditions \Regsoluza.
\Ethm

\Bdim
In order to prove the existence of a solution, 
we start from the finite-dimensional problem \Pbln, perform an a~priori estimate,
and let $n$ tend to infinity.
Also in this section, we simplify the notation as far as constants are concerned
and use the same symbol $c$ for different constants
that can depend only on the structure, the data, $T$, the optimal control $\ub$,
and the corresponding state $\soluzb$.

\step
A priori estimate

We write the equations \accorpa{priman}{terzan} at the time $s$ and test them by
$\,-\dt\qn(s)$, $\pn(s)$, and~\juerg{$\,-\dt\rn(s)$}, respectively.
Then, we sum up, integrate over~$(t,T)$ with respect to~$s$,
and notice that the terms involving the product $\pn\dt\qn$ cancel each other.
Moreover, we add the same quantities 
$\,\itT\norma\pn^2 \ds\,$ and $\,(1/2)\norma{\rn(t)}^2=\itT(\rn,\dt\rn)\ds\,$ to both sides 
in order to recover the full norms in the spaces $\VB\sigma$ and~$\VC\tau$.
We then obtain the identity
\begin{align}
  & \frac 1n \itT \norma{\dt\qn}^2 \ds
  + \frac 12 \, \norma{A^\rho\qn(t)}^2 
  + \frac 12 \, \norma{\pn(t)}^2
  + \itT \norma\pn_\Bs^2 \ds
  \non
  \\
  & \quad {}
  + \itT \norma{\dt\rn}^2 \ds
  + \frac 12 \, \norma{\rn(t)}_\Ct^2
  \non
  \\
  & {} = \itT \bigl( P(\phib)(\qn-\rn) \dt(\qn-\rn) \bigr) \ds
  + \itT \bigl( P'(\phib)(\Sb-\mub)(\qn-\rn) , \pn \bigr) \ds
  \non
  \\
  & \quad {}
  - \itT \bigl( f'(\phib)\pn , \pn \bigr) \ds
  + \itT (g_1,\pn) \ds
  - \itT (g_3,\dt\rn) \ds 
  \non
  \\
  & \quad {}
  + \frac 12 \, \norma{A^\rho\qn(T)}^2 
  + \frac 12 \, \norma{\pn(T)}^2
  + \frac 12 \, \norma{\rn(T)}_\Ct^2
  \non
  \\
  & \quad {}
  + \itT \norma\pn^2 \ds
  + \itT ( \rn,\dt\rn ) \ds \,.
  \label{perstimaFG}
\end{align}
\juerg{At first, we exploit the endpoint conditions \eqref{cauchyn}. 
Obviously, $\qn(T)=0$, which entails that $A^\rho\qn(T)=0$, as well as $(\pn(T),v)=(g_2,v)$ for all $v\in V_B^{(n)}$. 
The latter identity just means that
$\pn(T)$ is the $H$-orthogonal projection of $g_2$ onto $V_B^{(n)}$, which implies that $\,\|\pn(T)\|\,\le\,\|g_2\|\,$
for all $n\in\enne$. By the same token, we can infer that $\,\|\rn(T)\|\,\le\,\|g_4\|\,$ for all $n\in\enne$.
Finally, we insert $\,v=C^{2\tau}\rn(T)\in V_C^{(n)}\,$ in the last identity in \eqref{cauchyn}. 
Recalling that
$\overline{S}\in C^0([0,T];V_C^\tau)$, and by virtue of \eqref{compat}, we infer that $g_4\in V_C^\tau$. 
We thus find that}
\begin{align*}
&\juerg{\bigl\|C^\tau\rn(T)\bigr\|^2\,=\,\bigl(\rn(T),C^{2\tau}\rn(T)\bigr)\,=\,
\bigl(g_4,C^{2\tau}\rn(T)\bigr)\,=\,\bigl(C^\tau g_4,C^\tau\rn(T)\bigr),}
\end{align*}
\juerg{whence we infer that \,$\|C^\tau\rn(T)\|\,\le\,\|C^\tau g_4\|$. 
In conclusion, we have shown the estimate
$\|\rn(T)\|_{C,\tau}\,\le\,\|g_4\|_{C,\tau}\,$ for all $n\in\enne$.}

\juerg{Next, we consider the first two terms on the \rhs, which we denote by $Y_1$ and $Y_2$. We only need 
to estimate these terms, since the remaining other ones can easily be handled using Young's inequality and,
eventually, Gronwall's lemma.}
As for~$Y_1$, we first integrate by parts, 
and one of the important terms we obtain is nonpositive.
Then, we account for the \Holder\ and Youngs inequalities, 
the equivalence of norms in $\VA\rho$ related to~\eqref{forequiv},
and the embeddings~\eqref{embeddings} as follows:
\begin{align}
   Y_1\,& =\, \frac 12 \bintQt P(\phib) \, \dt|\qn-\rn|^2 \dx\ds
  \,=\, \frac 12 \iO P(\phib(T)) \, |\qn(T)-\rn(T)|^2
  \non
  \\
  & \quad  -\, \frac 12 \iO P(\phib(t)) \, |\qn(t)-\rn(t)|^2
  - \bintQt P'(\phib) \, \dt\phib \, |\qn-\rn|^2
  \non
  \\
  & \leq \,c + c \itT \norma{\dt\phib}_4 \, ( \norma\qn_4 + \norma\rn_4 )^2 \ds
  \,\leq \,c + c \itT \bigl( \norma{A^\rho\qn}^2 + \norma\rn_\Ct^2 \bigr) \ds \,.
  \non
\end{align}
Concerning $Y_2$, we have that 
\begin{align}
  Y_2 \,&\leq \,c \itT \norma{\Sb-\mub}_4 \, \norma{\qn-\rn}_4 \, \norma\pn_4 \ds
  \non
  \\ 
  &\leq\, \frac 12 \itT \norma\pn_\Bs^2 \ds
  + c \itT \bigl( \norma{A^\rho\qn}^2 + \norma\rn_\Ct^2 \bigr) \ds \,.
  \non
\end{align}
By treating the remaing terms on the \rhs\ of \eqref{perstimaFG} 
as announced before, and applying the Gronwall lemma, we conclude that
\begin{align}
  & \frac 1{\sqrt n} \, \norma{\dt\qn}_{\L2H}^2
  + \norma\qn_{\L\infty{\VA\rho}}
  \non
  \\[1mm]
  &   + \norma\pn_{\L\infty H\cap\L2{\VB\sigma}}
  + \norma\rn_{\H1H\cap\L\infty{\VC\tau}}
  \,\leq c \,.
  \label{stimaFG}
\end{align}

\step 
Existence

The above estimate ensures that, for a subsequence again indexed by $n$,
\begin{align}
  & \frac 1n \, \dt\qn \to 0 
  \quad \hbox{strongly in $\L2H$},
  \label{convdtqn}
  \\
  & \qn \to q 
  \quad \hbox{weakly star in $\L\infty{\VA\rho}$},
  \label{convqn}
  \\
  & \pn \to p
  \quad \hbox{weakly star in $\L\infty H\cap\L2{\VB\sigma}$},
  \label{convpn}
  \\
  & \rn \to r 
  \quad \hbox{weakly star in $\H1H\cap\L\infty{\VC\tau}$}.
  \label{convrn}
\end{align}
We aim at proving that $\soluza$ is the desired solution to the weak form \Pblav\ of the adjoint problem.
Clearly, \eqref{cauchyav}~is satisfied, and we have to prove that the variational equations are satisfied as well.
We confine ourselves to the second equation, which is the most complicated one.
To this end, we write an integrated version of~\eqref{secondan}.
We fix any integer $m>1$, take any $v\in\H1H\cap\L2\VBm$ vanishing at $t=0$, and assume that $n\geq m$.
Then\, $\VAm\subset\VAn$, so that $v(s)$ is admissible in~\eqref{secondan} written at the time~$s$,
and we can test the equation in the inner product of~$H$.
Moreover, we replace $(B^{2\sigma}\pn(s),v(s))$ by $(B^\sigma\pn(s),B^\sigma v(s))$.
Then, we integrate over $(0,T)$ with respect to~$s$.
Now, we observe that $v(T)$ is admissible in the second identity of~\eqref{cauchyn}.
So, by an integration by parts, we obtain that
\begin{align}
  & \ioT \bigl\{
    (\qn+\pn , \dt v) 
    + ( B^\sigma \pn , B^\sigma v ) + \bigl( f'(\phib) \, \pn - P'(\phib) (\Sb-\mub) (\qn-\rn) , v \bigr)
  \bigr\} \ds
  \non
  \\
  & {}
  = \ioT (g_1,v) \ds + \bigl( g_2 , v(T) \bigr)\,.
  \non
\end{align}
Since $n\geq m$ is arbitrary, we can let $n$ tend to infinity by using \accorpa{convqn}{convrn}.
Concerning, e.g., the worst term,
we recall that $\Sb$ and $\mub$ belong to $\L\infty{\Lx4}$ and observe that
$\qn$ and $\rn$ converge to $q$ and~$r$, respectively, also weakly in $\L2{\Lx4}$.
Hence, we have that
\Beq
  P'(\phib) (\Sb-\mub) (\qn-\rn) \to P'(\phib) (\Sb-\mub) (q-r)
  \quad \hbox{weakly in $\L2H$}.
  \non
\Eeq
As the other terms are easier, we conclude that \eqref{secondaav} holds for such a function~$v$.
At this point, we fix any $v\in\H1H\cap\L2{\VB\sigma}$ vanishing at $t=0$ and define $v_m$ by setting
\Beq
  v_m(t) := \somma j1m (v(t),e'_j) e'_j
  \quad \hbox{for $t\in[0,T]$}.
  \non
\Eeq
Then, $v_m$ belongs to $\H1H\cap\L2\VAm$ and vanishes at $t=0$.
Hence, we can use it in the equality just obtained.
As $m$ is arbitrary, we can \juerg{take the limit as $m\to\infty$}.
By noting that $v_m$ converges even strongly to $v$ in $\H1H\cap\L2{\VB\sigma}$,
we conclude that \eqref{secondaav} is satisfied for such a~$v$.
By similarly reasoning for the other equations, we can conclude.
Hence, the existence part of the statement is proved.

\step
Uniqueness

By linearity, we can assume that all the \rhs s of the strong formulation \Pbla\ vanish,
so that the problem becomes
\begin{align}
  & A^{2\rho} q - p + P(\phib) (q-r) 
  = 0 \,,
  \label{primau}
  \\[0.5mm]
  & - \dt(q+p) + B^{2\sigma} p  
  + f'(\phib) \, p - P'(\phib) (\Sb-\mub) (q-r)
  = 0\,,
  \qquad
  \label{secondau}
  \\[0.5mm]
  & - \dt r + C^{2\tau} r - P(\phib) (q-r) 
  = 0\,,
  \label{terzau}
  \\[0.5mm]
  & (q+p)(T) = 0
  \aand 
  r(T) = 0 \,.
  \label{cauchyu}
\end{align}
We cannot adapt the argument used to arrive at~\eqref{stimaFG},
since no information for $\dt q$ is available now.
So, we proceed in a different way.
With the notation
\Beq
  (1*v)(t) := \itT v(s) \ds
  \quad \hbox{\aat\ and every $v\in\L1H$}\,,
  \non
\Eeq
we integrate (\aat) \eqref{secondau} over $(t,T)$ and obtain \aet
\Beq
  q + p + B^{2\sigma}(1*p)
  = 1 * \bigl( P'(\phib) (\Sb-\mub) (q-r) \bigr)
  - 1 * \bigl( f'(\phib) \, p \bigr) \,.
  \label{intsecondau}
\Eeq
At this point, we test \eqref{primau} by~$q$, \eqref{intsecondau} by~$p$, and \eqref{terzau} by~$r$,
sum up, and integrate over~$(t,T)$.
The terms involving the product $p\,q$ cancel each other.
We also add the same quantities $\,(1/2)\norma{(1*p)(t)}^2=\itT p(1*p)\ds\,$ and $\,\itT\norma r^2\ds\,$
to both sides and obtain
\begin{align}
  & \itT \norma{A^\rho q}^2 \ds
  + \itT \bigl( P(\phib) (q-r) , q-r \bigr) \ds
  + \itT \norma p^2 \ds
  + \frac 12 \, \norma{(1*p)(t)}_\Bs^2
  \non
  \\
  & \quad {}
  + \frac 12 \, \norma{r(t)}^2
  + \itT \norma r_\Ct^2 \ds
  \non
  \\
  & {} = \itT \bigl(
    1 * \bigl( P'(\phib) (\Sb-\mub) (q-r) - f'(\phib) \, p \bigr) , p
  \bigr) \ds 
  \non
  \\
  & \quad {}
  + \itT p (1*p) \ds
  + \itT \norma r^2 \ds \,.
  \label{perunicita}
\end{align}
All of 	the terms on the \lhs\ are nonnegative.
Now, we treat the first integral on the \rhs, which we term~$Y$.
We first integrate by parts.
Then, we owe to Young's inequality and to the obvious inequality
$ \,\,\norma{(1*v)(t)}^2\leq T\itT\norma{v(s)}^2\ds$,
which holds true for every $t\in[0,T]$ and every $v\in\L2H$.
We set, for brevity, $w:=P'(\phib)(\Sb-\mub)(q-r)-f'(\phib)p$ \,and \pier{observe~that
\begin{align}
  Y\,& \leq \frac 1 4  \itT \norma p^2 \ds +  \itT \norma{1*w}^2 \ds  \non \\
   & \leq\,  \frac 1 4  \itT \norma p^2 \ds +  \itT T \int_s^T  \norma{w}^2 ds' \/ \ds.
  \non
\end{align}
On} the other hand, we recall \eqref{embeddings}, \eqref{forequiv} and the regularity $\L\infty{\Lx4}$ of $\Sb$ and~$\mub$.
We thus deduce \pier{that
\begin{align}
  &  \int_s^T  \norma{w}^2 ds' 
  \,\leq\,\norma{P'(\phib)}_\infty^2  \int_s^T \norma{\Sb-\mub}_4^2 \, \norma{q-r}_4^2 ds'
  + \norma{f'(\phib)}_\infty^2  \int_s^T \norma p^2 ds'
  \non
  \\
  & \leq\, c  \int_s^T \Bigl( \norma{A^\rho q}^2 + \norma r_\Ct^2 + \norma p^2 \Bigr) ds' \,,
  \non
\end{align}
whence
\Beq
  Y
  \leq 
   \frac 1 4  \itT \norma p^2 \ds +  \itT c \biggl( \int_s^T   \Bigl( \norma{A^\rho q}^2 + \norma p^2 + \norma r_\Ct^2  \Bigr)ds' \biggr) \ds.
\Eeq
Therefore,} coming back to \eqref{perunicita} and estimating the second integral on the \rhs\ \pier{as
\begin{align}
     \itT p (1*p) \ds \leq \frac 1 4  \itT \norma p^2 \ds +  \itT \norma{1*p}^2 \ds  
  \non
\end{align}
and then applying the Gronwall lemma,  we easily conclude that} $(q,p,r)=(0,0,0)$.
\Edim

\juerg{Now that} the adjoint problem is solved, we can rewrite the variational inequality \eqref{badnc}
in a much better form. \juerg{Indeed, we have the following result.}

\Bthm
\label{Goodnc}
Under the assumptions \HPcontrol\ \juerg{and \eqref{compat}}, let $\ub\in\Uad$ be an optimal control,
and let $\soluza$ be the solution to the associated adjoint problem \Pbla.
Then \juerg{it holds}
\Beq
  \intQ (r+\kappa_5\ub)(v-\ub) \geq 0
  \quad \hbox{for every $v\in\Uad$}.
  \label{goodnc}
\Eeq
In particular, if $\juerg{\kappa}_5>0$, then $\ub$ is the projection of $-r/\kappa_5$ on $\Uad$
in the sense of the space $\LQ2$ with its standard inner product. That is, it is given~by
\Beq
  \ub = \min \graffe{ \umax , \max \graffe{\umin,-r/\kappa_5} } \quad \aeQ \,.
  \non
\Eeq
\Ethm

\Bdim
We fix $v\in\Uad$ and consider the linearized system with $h=v-\ub$.
Now, we observe that the regularity \Regsoluza\ is suitable
for integrating over $(0,T)$ the equations \eqref{primal}, \eqref{secondal}, and \eqref{terzal},
tested by $q(t)$, $p(t)$, and~$r(t)$, respectively.
By doing this, rearranging and summing up, 
we~obtain (as~before in this section, we omit the integration variable, which we term $s$ for uniformity)
\begin{align}
  & \ioT \bigl\{
  \bigl( \dt\xi , q \bigr)
  + \bigl( A^\rho \eta , A^\rho q \bigr) 
  - \bigl( P(\phib) (\zeta-\eta) , q \bigr)
  - \bigl( P'(\phib) \, \xi \, (\Sb-\mub) , q \bigr)
  \bigr\} \ds
  \non
  \\
  & \quad {} + \ioT \bigl\{
  \bigl( \dt\xi , p \bigr)
  + \bigl( B^\sigma\xi , B^\sigma p \bigr)
  + \bigl( f'(\phib) \, \xi , p \bigr)
  - \bigl( \eta , p \bigr)
  \bigr\} \ds
  \non
  \\
  & \quad {} + \ioT \bigl\{
  \bigl( \dt\zeta , r \bigr)
  + \bigl( C^\tau \zeta , C^\tau r \bigr)
   + \bigl( P(\phib) (\zeta-\eta) , r \bigr)
  + \bigl( P'(\phib) \, \xi \, (\Sb-\mub) , r \bigr)
  \bigr\} \ds
  \non
  \\
  & = \ioT \bigl( v-\ub , r \bigr) \ds \,.
  \non
\end{align}
At the same time, we \pier{take $v = -\eta$ in \eqref{primaav}, $v=-\xi $ in \eqref{secondaav}, $v=-\zeta$ in \eqref{terzaav}, respectively, and note that all the three test functions are admissible} in their equations.
Then, we sum up, rearrange, and~get
\begin{align}
  & \ioT \bigl\{ - \bigl( A^\rho q , A^\rho \eta \bigr)
  + ( p,\eta ) - \bigl( P(\phib) (q-r) ,\eta ) \bigr\} \ds
  \non
  \\
  & \quad {}
  + \ioT \bigl\{
    - (q+p , \dt \xi) 
    - ( B^\sigma p , B^\sigma\xi ) - \bigl( f'(\phib) \, p - P'(\phib) (\Sb-\mub) (q-r) , \xi \bigr)
  \bigr\} \ds
  \qquad
  \non
  \\
  & \quad {}
  + \ioT \bigl\{
    ( \dt r , \zeta)
    - ( C^\tau r , C^\tau \zeta) + \bigl( P(\phib) (q-r) , \zeta \bigr)
  \bigr\} \ds
  \non
  \\
  & {}
  = - \ioT (g_1,\xi) \ds 
  - \bigl( g_2 , \xi(T) \bigr)
  - \ioT (g_3 , \zeta ) \ds \,.
  \non
\end{align}
Next, we add the identities just obtained to each other.
Several cancellations occur, and what remains is just the following identity:
\Beq
  \ioT \{ (\dt\zeta,r) + (\dt r,\zeta) \} \ds
  = \ioT \bigl( v-\ub , r \bigr) \ds 
  - \ioT (g_1,\xi) \ds 
  - \bigl( g_2 , \xi(T) \bigr)
  - \ioT (g_3 , \zeta ) \ds \,.
  \non
\Eeq
At this point, we observe that the \lhs\ equals $(g_4,\juerg{\zeta}(T))$ by~\eqref{cauchyav},
so that the above identity becomes
\Beq 
  \ioT (g_1,\xi) \ds 
  + \bigl( g_2 , \xi(T) \bigr)
  + \ioT (g_3 , \zeta ) \ds 
  + \bigl( g_4 , \zeta(T) \bigr) 
  = \ioT \bigl( v-\ub , r \bigr) \ds \,.
  \non
\Eeq
Hence, by recalling the notation \eqref{abbrev}, and comparing with~\eqref{badnc},
we obtain~\eqref{goodnc}.
\Edim

%%%%%%%%%%%%%%%%%%%%%%%%%%%%%%%%%%%%%%%%%%%%%%%%%%%%%%%%%%%%%%%%%%%%%%%%

\section*{Acknowledgments}
\pier{This research was supported by the Italian Ministry of Education, 
University and Research~(MIUR): Dipartimenti di Eccellenza Program (2018--2022) 
-- Dept.~of Mathematics ``F.~Casorati'', University of Pavia. 
In addition, PC and CG gratefully acknowledge some other 
financial support from the GNAMPA (Gruppo Nazionale per l'Analisi Matematica, 
la Probabilit\`a e le loro Applicazioni) of INdAM (Isti\-tuto 
Nazionale di Alta Matematica).}

%%%%%%%%%%%%%%%%%%%%%%%%%%%%%%%%%
%% bibliography
%%%%%%%%%%%%%%%%%%%%%%%%%%%%%%%%%

\vspace{3truemm}

\Begin{thebibliography}{10}

\bibitem{AM}
M. Ainsworth, Z. Mao,
Analysis and approximation of a fractional Cahn--Hilliard equation, \emph{SIAM J. Numer. Anal.} {\bf 55} (2017),
1689-1718.

\bibitem{ASS}
G. Akagi, G. Schimperna, A. Segatti,
Fractional Cahn--Hilliard, Allen--Cahnn, and porous medium equations,
\emph{J. Differential Equations} {\bf 261} (2016), 2935-2985.

\bibitem{Brezis}
H. Brezis,
``Op\'erateurs maximaux monotones et semi-groupes de contractions
dans les espaces de Hilbert'',
North-Holland Math. Stud.
{\bf 5},
North-Holland,
Amsterdam,
1973.

\bibitem{CRW}
C. Cavaterra, E. Rocca, H. Wu,
Long-time dynamics and optimal control of a diffuse interface model for tumor growth,
{\it Appl. Math. Optim.} (2019), https://doi.org/10.1007/s00245-019-09562-5.

\bibitem{CG}
P. Colli, G. Gilardi,
Well-posedness, regularity and asymptotic analyses for a fractional phase field system,
{\em Asymptot. Anal.}, to appear (see also preprint arXiv:1806.04625 [math.AP] (2018), pp. 1-34).

\bibitem{CGH}
P. Colli, G. Gilardi, D. Hilhorst,
On a Cahn--Hilliard type phase field system related to tumor growth,
{\it Discrete Contin. Dyn. Syst.} {\bf 35} (2015), 2423-2442.

\bibitem{SM}
P. Colli, G. Gilardi, G. Marinoschi, E. Rocca,
Sliding mode control for a phase field system related to tumor growth,
{\it Appl. Math. Optim.} \textbf{79} (2019), 647-670.

\bibitem{CGRS_VAN}
P. Colli, G. Gilardi, E. Rocca, J. Sprekels,
Vanishing viscosities and error estimate for a Cahn--Hilliard type phase field system related to tumor growth,
{\it Nonlinear Anal. Real World Appl.} {\bf 26} (2015), 93-108.

\bibitem{CGRS_OPT}
P. Colli, G. Gilardi, E. Rocca, J. Sprekels,
Optimal distributed control of a diffuse interface model of tumor growth,
{\it Nonlinearity} {\bf 30} (2017), 2518-2546.

\bibitem{CGRS_ASY}
P. Colli, G. Gilardi, E. Rocca, J. Sprekels,
Asymptotic analyses and error estimates for a Cahn--Hilliard type phase field system modelling tumor growth,
{\it Discrete Contin. Dyn. Syst. Ser. S} {\bf 10} (2017), 37-54.

\bibitem{CGS14}
P. Colli, G. Gilardi, J. Sprekels,
Optimal velocity control of a viscous Cahn--Hilliard system with convection and dynamic boundary conditions,
\juerg{\emph{SIAM J. Control Optim.}} 
{\bf 56} (2018), 1665-1691.

\pier{\bibitem{CGS18}
P. Colli, G. Gilardi, J. Sprekels,
Well-posedness and regularity for a generalized fractional Cahn--Hilliard system,
{\em Atti Accad. Naz. Lincei Rend. Lincei Mat. Appl.}, to appear 
(see also preprint arXiv:1804.11290 [math.AP] (2018), pp. 1-36).}

\bibitem{CGS19}
P. Colli, G. Gilardi, J. Sprekels,
Optimal distributed control of a generalized fractional Cahn--Hilliard system,
{\em Appl. Math. Optim.} \pier{(2019), %Online First 15 November  2018, 
https://doi.org/10.1007/s00245-018-9540-7.}

\bibitem{CGS21}
P. Colli, G. Gilardi, J. Sprekels,
Deep quench approximation and optimal control of general Cahn--Hilliard 
systems with fractional operators and double-obstacle potentials,
preprint arXiv:1812.01675 [math.AP] (2018), pp. 1-32.  

\bibitem{CGS23}
P. Colli, G. Gilardi, J. Sprekels,
Well-posedness and regularity for a fractional tumor growth model,
{\em Adv. Math. Sci. Appl.} {\bf 28} (2019), 343-375.

\bibitem{CGS22}
P. Colli, G. Gilardi, J. Sprekels,
Longtime behavior for a generalized Cahn--Hilliard system with fractional operators,
preprint arXiv:1904.00931 [math.AP] (2019), pp. 1-18.

\bibitem{CGS24}
P. Colli, G. Gilardi, J. Sprekels,
Asymptotic analysis of a tumor growth model with fractional operators,
in preparation.

\bibitem{ConGio}
M. Conti, A. Giorgini,
The three-dimensional Cahn--Hilliard--Brinkman system with unmatched densities,
preprint hal-01559179 (2018), pp. 1-34. 

\bibitem{CLLW}
V. Cristini, X. Li, J. S. Lowengrub, S. M. Wise,
Nonlinear simulations of solid tumor growth using a mixture model: invasion and branching.
{\it J. Math. Biol.} {\bf 58} (2009), 723-763.

\bibitem{CL}
V. Cristini, J. S. Lowengrub,
``Multiscale Modeling of Cancer: An Integrated Experimental and Mathematical
Modeling Approach'', Cambridge University Press, Leiden (2010).

\bibitem{DFRGM}
M. Dai, E. Feireisl, E. Rocca, G. Schimperna, M. Schonbek,
Analysis of a diffuse interface model of multi-species tumor growth,
{\it Nonlinearity\/} {\bf  30} (2017), 1639-1658.

\bibitem{DGG}
F. Della Porta, A. Giorgini, M. Grasselli,
The nonlocal Cahn--Hilliard--Hele--Shaw system with logarithmic potential,
{\em \pier{Nonlinearity}} {\bf 31} (2018), 4851-4881.

\pier{%
\bibitem{EGAR}
M. Ebenbeck, H. Garcke,
Analysis of a Cahn--Hilliard--Brinkman model for tumour growth with chemotaxis,
{\it J. Differential Equations,} \textbf{266} (2019), 5998-6036.
\bibitem{EK}
M. Ebenbeck, P. Knopf,
Optimal medication for tumors modeled by a Cahn--Hilliard--Brinkman equation,
preprint arXiv:1811.07783 [math.AP] (2018), pp. 1-26.
\bibitem{EK_ADV}
M. Ebenbeck, P. Knopf,
Optimal control theory and advanced optimality conditions 
for a diffuse interface model of tumor growth,
preprint arXiv:1903.00333 [math.OC] (2019), pp. 1-34.}

\bibitem{FGR}
S. Frigeri, M. Grasselli, E. Rocca,
On a diffuse interface model of tumor growth,
{\it  European J. Appl. Math.\/} {\bf 26} (2015), 215-243. 

\bibitem{FLRS}
S. Frigeri, K. F. Lam, E. Rocca, G. Schimperna,
On a multi-species Cahn--Hilliard--Darcy tumor growth model with singular potentials,
{\it Commun. Math Sci.} {\bf  (16)} (2018), 821-856. 

\bibitem{Gal1}
C. G. Gal, On the strong-to-strong interaction case for doubly nonlinear Cahn--Hilliard equations,
\emph{Discrete Contin. Dyn. Syst.} {\bf 37} (2017), 131-167.

\bibitem{Gal2}
C. G. Gal, Non-local Cahn--Hilliard equations with fractional dynamic boundary conditions,
\emph{European J. Appl. Math.} {\bf 28} (2017), 736-788.

\bibitem{Gal3}
C. G. Gal, Doubly nonlinear  Cahn--Hilliard equations, \emph{Ann. Inst. H. Poincar\'e Anal. Non Lin\'eaire}
{\bf 35} (2018), 357-392.

\bibitem{GARL_3}
H. Garcke, K. F. Lam,
Global weak solutions and asymptotic limits 
of a Cahn--Hilliard--Darcy system modelling tumour growth,
{\it AIMS Mathematics} {\bf 1} (2016), 318-360.
\bibitem{GARL_1}
H. Garcke, K. F. Lam,
Well-posedness of a Cahn--Hilliard system modelling tumour
growth with chemotaxis and active transport,
{\it European. J. Appl. Math.} {\bf 28} (2017), 284-316.
\bibitem{GARL_2}
H. Garcke, K. F. Lam,
Analysis of a Cahn--Hilliard system with non--zero Dirichlet 
conditions modeling tumor growth with chemotaxis,
{\it Discrete Contin. Dyn. Syst.} {\bf 37} (2017), 4277-4308.
\bibitem{GARL_4}
H. Garcke, K. F. Lam,
On a Cahn--Hilliard--Darcy system for tumour growth 
with solution dependent source terms, 
in ``Trends on Applications of Mathematics to Mechanics'', 
E.~Rocca, U.~Stefanelli, L.~Truskinovski, A.~Visintin~(ed.), 
{\it Springer INdAM Series} {\bf 27}, Springer, Cham, 2018, pp. 243-264.
\bibitem{GAR}
H. Garcke, K. F. Lam, R. N\"urnberg, E. Sitka,
A multiphase Cahn--Hilliard--Darcy model for tumour growth with necrosis,
{\it Math. Models Methods Appl. Sci.} {\bf 28} (2018), 525-577.

\bibitem{GARLR}
H. Garcke, K. F. Lam, E. Rocca,
Optimal control of treatment time in a diffuse interface model of tumor growth,
{\it Appl. Math. Optim.} {\bf 78} (2018), {495-544}.

\bibitem{GLSS}
H. Garcke, K. F. Lam, E. Sitka, V. Styles,
A Cahn--Hilliard--Darcy model for tumour growth with chemotaxis and active transport,
{\it Math. Models Methods Appl. Sci. } {\bf 26} (2016), 1095-1148.

\bibitem{HDPZO}
A. Hawkins-Daarud, S. Prudhomme, K. G. van der Zee, J. T. Oden,
Bayesian calibration, validation, and uncertainty quantification of diffuse 
interface models of tumor growth. 
{\it J. Math. Biol.} {\bf 67} (2013), 1457-1485. 

\bibitem{HDZO}
A. Hawkins-Daruud, K. G. van der Zee, J. T. Oden, Numerical simulation of
a thermodynamically consistent four-species tumor growth model, 
{\it Int. J. Numer. Math. Biomed. Engng.} {\bf 28} (2012), 3-24.

\bibitem{Lions}
J.-L. Lions, ``Quelques M\'ethodes de R\'esolution des Probl\`emes aux Limites non Lin\'eaires'', 
\pier{Dunod}, Gauthier-Villars, Paris, 1969. 

\bibitem{MRS}
A. Miranville, E. Rocca, G. Schimperna,
On the long time behavior of a tumor growth model,
{\it J. Differential Equations\/} {\bf 267} (2019), 2616-2642.

\bibitem{OHP}
J. T. Oden, A. Hawkins, S. Prudhomme,
General diffuse-interface theories and an approach to predictive tumor growth modeling,
{\it Math. Models Methods Appl. Sci.} {\bf 20} (2010), 477-517. 

\pier{%
\bibitem{S}
A. Signori,
Optimal distributed control of an extended model of tumor 
growth with logarithmic potential,
{\it Appl. Math. Optim.} (2018), https://doi.org/10.1007/s00245-018-9538-1.
\bibitem{S_DQ}
A. Signori,
Optimality conditions for an extended tumor growth model with 
double obstacle potential via deep quench approach, 
preprint arXiv:1811.08626 [math.AP] (2018), pp. 1-25.
\bibitem{S_b}
A. Signori,
Optimal treatment for a phase field system of Cahn--Hilliard 
type modeling tumor growth by asymptotic scheme, 
preprint arXiv:1902.01079 [math.AP] (2019), pp. 1-28.
\bibitem{S_a}
A. Signori,
Vanishing parameter for an optimal control problem modeling tumor growth,
preprint arXiv:1903.04930 [math.AP] (2019), pp. 1-22.%
}

\bibitem{SW}
J. Sprekels, H. Wu,
Optimal distributed control of a Cahn--Hilliard--Darcy system with mass sources,
{\it Appl. Math. Optim.} (2019), https://doi.org/10.1007/s00245-019-09555-4.

\bibitem{WLFC}
S. M. Wise, J. S. Lowengrub, H. B. Frieboes, V. Cristini,
Three-dimensional multispecies nonlinear tumor growth - I: model and numerical method, 
{\it J. Theor. Biol.} \pier{{\bf 253}} (2008), 524-543.

\bibitem{WZZ}
X. Wu, G. J. van Zwieten, K. G. van der Zee, Stabilized second-order splitting
schemes for Cahn--Hilliard models with applications to diffuse-interface tumor-growth models, 
{\it Int. J. Numer. Meth. Biomed. Engng.} {\bf 30} (2014), 180-203.

\End{thebibliography}

\End{document}

%%%%%%%%%%%%%%%%%%%%%%%%%%%%%%%%%%%%%%%%%%%%%